\documentclass[12pt,a4paper]{article}
 
 \usepackage{enumerate}
 \usepackage{subfigure}
 \usepackage{printlen}
 \usepackage{graphicx}

\usepackage{amsfonts,amssymb,amsmath} 

\usepackage{hyperref}
\usepackage{changes}
\usepackage{a4wide}
\usepackage{anysize}
\newtheorem{theorem}{Theorem}
\usepackage{times}
\usepackage{multirow}
\usepackage{wrapfig}
\usepackage{fullpage}
\usepackage{rotating}
\usepackage{tikz}
\usepackage{hyperref}
\usepackage{changes}
\usepackage{enumerate}
\newtheorem{result}{Result}

\newtheorem{lemma}{Lemma}
\newtheorem{conjecture}{\bf Conjecture}

\usepackage{euscript}
 \newcommand{\eop}{{\hfill $\blacksquare$} }
 
 \newcommand{\lesseq}{ \hspace{-1.5mm}&  \hspace{-1.5mm} :=  \hspace{-1.5mm}& \hspace{-1.5mm}}
 
 \newcommand{\kav}[1]{\textcolor{red}{#1}}
 
 \newcommand{\Kmin}{{\underline k}^s}

 \newcommand{\Kmax}{{\overline k}^s}

 \usepackage{ifthen}

\newcommand{\ForTR}{\ifthenelse{1<2}} 
\newcommand{\ATR}{\ifthenelse{1<2}}

\newcommand{\x}{{\bar {\bf x}}}
\newcommand{\X}{{\bar {\bf X}}}

\newcommand{\kz}{k_z}
\newcommand{\y}{{\tilde y}^b}

\renewcommand{\u}{{\bar {\bf u}}} 

\newcommand{\ignore}[1]{}

\ForTR{
\newcommand{\PhiPsi}{\Psi}
}{
\newcommand{\PhiPsi}{\Phi}}

\newcommand{\f}{{\bf  f}}
\begin{document}

%\frontmatter
\thispagestyle{empty}
 \title{Random Fixed Points,  Limits and  
 %Heterogeneous 
 Systemic risk }

 \author{$^1$Veeraruna Kavitha, $^1$Indrajit Saha and  $^2$Sandeep Juneja \\  
 $^1$IEOR, IIT Bombay,  and 
 $^2$TIFR Mumbai,   India}

 \maketitle

\begin{abstract}
We consider vector fixed point (FP) equations in large dimensional spaces involving random variables, and study their  realization-wise solutions. 
We have an underlying directed random graph, that defines the connections between various components of the FP equations.  Existence of an edge between  nodes  $i,j$ implies the $i$-th FP 
equation  depends on the $j$-th component. 
We consider a special case where any component of the FP equation depends upon an appropriate aggregate of that of the random `neighbour'  components. We obtain   finite dimensional limit FP equations (in a much   smaller dimensional space), whose solutions approximate the solution of the random FP equations for   almost all realizations,  in the asymptotic limit (number of components increase).  
Our techniques are  different from the traditional mean-field  methods,  which  deal with
 stochastic FP equations in the space of distributions   to describe  the   stationary distributions of the systems. In contrast  our focus is on realization-wise FP solutions.  
We apply the results to study systemic risk in a large financial \underline{heterogeneous} network with many small institutions and one big institution, and demonstrate some interesting phenomenon.

\end{abstract}

 \section{Introduction}
 
 Random fixed points (FPs) are generalization of classical deterministic  FPs,  and arise  when one considers systems with uncertainty.  
 Broadly one can consider two types of such fixed points. There is considerable literature that  considers  stochastic FP equations on the space of probability distributions 
 (e.g., \cite{Urns,WeightedBranch}).   These equations typically arise as a limit of some iterative schemes, or as asymptotic (stationary) distribution of stochastic systems.  Alternatively, one might be interested in sample path wise FPs (e.g., \cite{Measure_FP,Measure_Approx}).  For each realization of the random quantities describing the system, we have one deterministic FP equation. These kind of equations can arise  when the performance/status of an  agent  depends upon that of a number of other agents. For example, a financial network with any given liability graph   is affected by individual/common random economic shocks received by the agents.   
 The amount cleared (full/fraction of liability) by an agent  
 depends upon: a) the   shocks  it receives; and  b) the liabilities cleared by the other agents. %The infection status  of an individual in a polluted environment depends upon the individuals immunity levels as well as the infection status of its neighbours.   
 Our primary focus in this paper is on the second type of equations.  Existing literature primarily considers the existence  of measurable FP, given the existence of realization-wise   FPs (e.g., \cite{Measure_FP,Measure_Approx}). In \cite{Measure_Approx} (and reference therein) authors consider the idea of random proximity points.

 To the best of our knowledge, there are no (common) techniques that provide `good' solutions to (some special types of) these equations. 
 We consider a special type of FP equations,  which are quite common, and provide a procedure to compute the approximate solutions. 
 We have FP equations
 in which the performance/status of an agent is influenced only by the aggregate  performance/status of its neighbours. 
 A random graph describes the neighbours, while a set of FP equations (one per realization of the random quantities) describe the  performance vectors.  
The key  idea is to study these FPs, asymptotically as the number of agents increase.    Towards this, we first study the aggregate influence factors, 
with an aim to reduce the  dimensionality of the problem. But due to random connections,  
the aggregate influence  factors can also depend upon the nodes.  
  However   the  aggregates might converge towards the same limit (e.g., as in  law of large numbers). % , one can anticipate a fewer  distinct aggregate factors at limit.  
 {We  precisely consider such scenarios and show that the random FPs}   converge to  that of  a limit system, under certain conditions.
\ForTR{ The performance of  the agents   in the limit system, depends upon finitely many  `aggregate'  limits. }{}We could also obtain   closed form expressions for   approximate  solutions of some   examples.

  \ForTR{The mean field theory (MFT) is close to this approach: 
MFT approximates many body problem with a one body problem and our result is also similar in nature.  However there are significant differences.}{}  
The mean field theory    also deals with a system of large number of agents, wherein the state/behaviour of an individual agent  is influenced by its own (previous) state and the mean (aggregate)  field  seen by it (e.g., \cite{Mean_WLAN} and reference therein).\ForTR{ The mean field is largely described in terms of  occupation (empirical) measures representing the fraction of agents in different states.}{}  
 The theory shows the convergence of the state trajectories as well as the stationary (time limit) distributions of the  original system towards that of a limit deterministic system. 
%The state of an individual agent, in the limit system, evolves depending upon its own state and the deterministic trajectories obtained by the limit of the  
%mean field trajectories. 
 %However our system is directly described by a set of FP equations, which depend upon the `mean' performance. 
  The stationary distribution   can be described by  FP equations in the space of distributions (e.g., \cite{Mean_WLAN}).  While 
we directly have a set of FP equations, which are defined realization-wise and depend upon the realization-wise \ `mean' performance.   
  Further, as already mentioned the mean influence factor is not common to all the agents.  In \cite{Mean_rand} and references therein, authors consider the mean field analysis  with `random' aggregate influence factors like in our case. 
  They consider the first-order   approximation, wherein the   joint expected values are approximated
by the product of the  marginal expected values etc.  Thus the moments of the joint distributions representing the FP solutions are asymptotically proven to be product of marginal moments.     Some authors also consider second-order approximations or moment closure techniques, where
the joint states of triplets are assumed to have a specific distribution (see \cite{Mean_rand} for relevant discussions). Our FP solutions are also proved to be asymptotically independent, however the asymptotic solutions are  
independent (infinite dimensional) random vectors.
 
We consider FP equations with possibly multiple solutions. We show that any sequence of the chosen FPs, converges   to the unique FP of the limit system almost surely\ForTR{ under sufficiently general conditions.}{.}  
%Secondly we have random FPs, and one could have used the above approach if we were interested in  only sample path-wise results (and in case they were strict contractions).   But 
   
 We   apply our results to study the systemic risk  in a  large financial network with many financial institutions. 
The institutions borrow/lend money from/to other institutions, and will have to clear their obligations at a later time point.  These systems are subjected to economic shocks, 
because of which some entities default (do not clear their obligations).  Because of interdependencies, this can lead to further defaults and the cascade of these reactions can lead to  
(partial/full)  {\it collapse} of the system.  %Systemic risk deals precisely with 
% One can define the systemic risk as the risk of 'collapse' of  the system because of enconmic shocks, as a result of the actions taken by  individual/group of components.  
 After the financial crisis of $2007$-$2008$, there is a surge of activity towards studying systemic risk (e.g.,  \cite{allen2000financial},\cite{acemoglu2015systemic},\cite{eisenberg2001systemic}).  The focus in these papers has been on several aspects including, measures to capture systemic risk,   influence of network structure on systemic risk, phase transitions etc.  These papers primarily discuss homogeneous systems, although  heterogeneity 
 is a crucial feature of real world networks.    As already mentioned 
 the clearing vectors are represented by FP equations and may have multiple FP solutions.
Thus our asymptotic solution can be  useful in this context.    We consider one  stylized example of  \underline{heterogeneous financial } network, that of one big bank and numerous small banks. 
Our key contribution is that we develop a methodology to arrive at simplified asymptotic representation to large bank networks.
This allows easy resolution of many practical what-if  scenarios. For instance, in a simple framework we observe that
having a big bank in an economy well connected to the small banks can stabilize the small banks even when the big bank itself faces shocks. However,
the reverse may not be true.  
The proposed methodology can be similarly used to provide insights into many other practical scenarios. We  analyze these in future. To summarise,  our analysis helps identify important patterns in a complex structure, since the structure simplifies when large number of constituents are involved.
  
 % \vspace{-3mm}
%  In a technical report \cite{TR},  towards the end we briefly describe another application, that of approximately reducing the state space of a Markov decision process, under certain nearness assumptions.  The second one is a very brief discussion, only meant for illustration purpose of other applications of our results, and requires further elaboration. 
% 
% 

 \vspace{-3mm}
\section{System Model  }
%Consider an IID sequence of random variables/vectors $G_i$. The various  components of any $G_i$ can be correlated but $G_i$ for different $i $ are independent of each other. Consider another IID sequence of double indexed sequence $\{ I_{i,j} \}_{i,j}$ of weights,  where  $ I_{i,j} > 0$ indicates the existence of an edge between nodes $i, j$ in a sequence of graphs as discussed below. 

Consider a  random graph with $n+1$ vertices $\{1, 2, \cdots, n, b\}$  whose directed   edges,   given  by random weights    $\{W_{i,j}\}$, represent the influence factors.
The node $b$ is a `big node',   and is highly influential. 
There is an edge between any two of the `small' nodes (nodes in $\{1, 2, \cdots, n\}$) with probability $p_{ss}$ independently of the others 
and let 
 $\{ I_{i,j} \}_{i \le n, j\le n}$ be the corresponding indicators.   Then the weights from a small node $j$ are the fractions\footnote{Note that
 $\sum_{i} W_{j,i} + W_{j, b} = 1$ for all $j$.} defined as below:
 \vspace{-2.5mm}
\begin{eqnarray}
 \label{Eqn_Weights} \hspace{20mm}
W_{j, b} =  \eta_j^{sb} \mbox{ and }  W_{ j, i}  =  \frac{I_{j,i} (1-\eta_j^{sb}) } { \sum_{i' \le n} I_{j, i'}  } , 
\end{eqnarray}where $\{\eta_j^{sb} \}_j$ are IID (independent, identically distributed) random variables with values between $0,1$. These fractions, for example,  can represent random fractions of some resources shared between various nodes. From small node $j$, there is a dedicated fraction $\eta^{sb}_j$ towards the b-node while the remaining  ($1-\eta^{sb}_j$ ) fraction is equally shared by the other connected small nodes.    The weights from $b$-node are the fractions, 
$$
\hspace{19mm}
W_{b, j}  = \frac{\eta_j^{bs} }{ \sum_i \eta_i^{bs}},
$$ where $\{\eta_j^{bs} \}_j$ are IID random variables again.  We are interested in   some  performance of the nodes, which depends upon the weighted average of the performance  of   other nodes with weights as given by $\{W_{i, j}\}.$
We consider the following fixed point (FP) equation  (in $R^{n+1}$)  constructed using  functions $(f^s, f^b)$, which in turn depend 
upon weighted averages $\{{\bar X}^s_i \}_i$ and ${\bar X}^b$,  and 
whose FP  ($i$-th component)
 represents   important performance measure of the nodes (node-$i$) as below:

{\small \begin{eqnarray}
\label{Eqn_FixedeqGen}
X_i^{s}  &=& f^s(G_i,   {\bar X}^s_i , \eta_i^{bs} X^b) \mbox{ {\normalsize  for each }} i \le n,  \\   
\label{Eqn_Fixedeq2Gen}X^b  &=&    f^b (   {\bar X}^b ) \ \ \ \mbox{ {\normalsize with aggregates}}    \\
{\bar X}^s_i  & := &  \sum_{j \le n} X_j^s W_{j, i}  \  \mbox{ {\normalsize and} } \  \nonumber 
{\bar X}^b :=  \frac{1}{n}  \sum_{j \le n} {X}^s_j W_{j, b} .
\end{eqnarray}}%
In the above, $\{G_i\}$ is an IID sequence and the performance of the big node $X^b$ is defined  per small node (performance divided by $n$).  
For any $n$ define mapping $\f := (f^b, f^s,   \cdots f^s)$, with ${\bf x}:= {\bf x}^n := (x_1^n, x_2^n, \cdots, x_n^n)$,    component wise:   

\vspace{-8mm}
{\small \begin{eqnarray*}
f_{1}( {\bf x}, x_b ) \hspace{-1mm}&\hspace{-1mm} :=\hspace{-1mm}& \hspace{-1mm} f^b ( {\bar x}_b),  \   {\bar x}_b  :=  \frac{1}{n}\sum_{j \le n} x_j     W_{j,b} \mbox{ and} \\
f_i ( {\bf x}, x_b ) \hspace{-1mm}&\hspace{-1mm} :=\hspace{-1mm}& \hspace{-1mm}   f^s(G_i,  {\bar x}_i, \eta_i^{bs}x_b ) , \   {\bar x}_i  :=  \sum_{j \le n} x_j     W_{j,i} 
\ \forall  i > 1,
\vspace{-1mm}
\end{eqnarray*}}which represents the FP of the random  operator (\ref{Eqn_FixedeqGen})-(\ref{Eqn_Fixedeq2Gen}). 
We assume the following:\\
{\bf A.1}
  The functions $f^s, f^b$ are non-negative, continuous and are bounded by an $y < \infty$, % i.e., for e.g.,
  
{\small $$
0 \le f^s(g,    x, x_b), \ f^b (x_b) \le y \mbox{ for all }  g, x, x_b. 
$$}
 Under the above assumption, by well known Brouwers fixed point theorem, FP solution exists for     almost all realizations of $\{G_i\}$,  $\{W_{j,i}\}$ and for any $x_b$.  Thus we have a random  (measurable) FP $( {\bf X}^{*}, X_b^*)$ for each $n$ for the random  operator (\ref{Eqn_FixedeqGen})-(\ref{Eqn_Fixedeq2Gen}) (see \cite{Measure_FP}). \ATR{To be precise we have:}{The precise details are as Lemma 1 of \cite{TR}.}
\ATR {
\begin{lemma}
\label{CLemma_exist}
For any $n$ define mapping $\f := (f^b, f^s,   \cdots f^s)$, with ${\bf x}:= {\bf x}^n := (x_1^n, x_2^n, \cdots, x_n^n)$,    component wise:   

\vspace{-9mm}
{\small \begin{eqnarray*}
f_{1}( {\bf x}, x_b ) \hspace{-1mm}&\hspace{-1mm} :=\hspace{-1mm}& \hspace{-1mm} f^b ( {\bar x}_b),  \   {\bar x}_b  :=  \frac{1}{n}\sum_{j \le n} x_j     W_{j,b} \mbox{ and} \\
f_i ( {\bf x}, x_b ) \hspace{-1mm}&\hspace{-1mm} :=\hspace{-1mm}& \hspace{-1mm}   f^s(G_i,  {\bar x}_i, \eta_i^{bs}x_b ) , \   {\bar x}_i  :=  \sum_{j \le n} x_j     W_{j,i} 
\ \forall  i > 1.
\vspace{-1mm}
\end{eqnarray*}}
Each component 
 is a mapping from $[0, y ]^{n+1}  \to [0, y ]$   
  for almost all $\{G_i\}$,  $\{W_{j,i}\}$ and for any $x_b$.   
 Further by continuity of $\f$,  we have a deterministic fixed point   for almost all $\{G_i\}$   and $\{I_{j,i}\}$ under {\bf A}.1.  Then we have (almost sure) random fixed point $( {\bf X}^{*}, X_b^*)$ for each $n$ (see \cite{Measure_FP}).  \eop
 \end{lemma}
}{}

{\bf Assumptions on the graph structure:}  We require that the number of  nodes  influencing any given  node,  grows asymptotically linearly  for almost all   sample paths:
\noindent {\bf A.2}   Consider $p_{ss} > 0$, and  only graphs for which,
 %
 %\vspace{-8mm}
 {\small
\begin{eqnarray*}
  \lim_{n \to \infty }  \sum_{j \le n}    \left|  
    \frac{1 }{   \sum_{i} I_{j, i}}  - \frac{1}{ n p_{ss}}    \right |  = 0   \mbox{ almost surely (a.s.),       for any } i.     
\end{eqnarray*}}

\vspace{-5mm}
\subsection{ Aggregate fixed points}
One can rewrite the fixed point equations for the weighted averages  $\{ {\bar X}_i^s \}_i$, ${\bar X}^b$ and we begin with their analysis. 
 Define the following random variables, that depend upon real constants $(x, x_b)$:\vspace{-5mm}
\begin{eqnarray}\hspace{15mm} \label{Eqn_xi}
\xi_i (x, x_b) :=   f^s (G_i,   x , \eta_{i}^{bs} x_b) ,
 %   \mbox{ with } l := 1 - E[\eta^{sb}_i] ,
\end{eqnarray} and assume: \\
 {\bf A.3}   {\small $|\xi_i (x, x_b)  - \xi_i (u, u_b) | \le \sigma ( |x- u| + |x_b- u_b|)$} with  {\small $\sigma \le 1$}. 
 
Consider  the  following  operators   on infinite  sequence space\footnote{
\label{footnote_sinf}
 Here $s^\infty$ is the space (subset) of  bounded sequences equipped with $l^\infty$ norm  $ |\x |_\infty :=  \sup_i |x_i|$, 
$$
s^\infty := \{  \x =  (x_1, x_2, \cdots ) : x_i \in [0, y ]  \mbox{ for all }  i    \}.
$$} $s^\infty$,   one for each $ n$: 

\vspace{-8mm}
{\small \begin{eqnarray}
\label{Eqn_bar_fixedpoint_randGen}
{\bar {\bf f}}^n ({\x}, {\bar x}_b)  &=&  ({\bar f}^n_b, {\bar f}^n_1, {\bar f}^n_2 \cdots) \mbox{ where} \\ \nonumber
{\bar f}^n_i ( \x, {\bar x}_b ) &=&  
\left \{  \begin{array}{lll}
 \sum_{j \le n}   \xi_{j} ({\bar x}_j, x_b)  W_{j, i}  &\mbox{ if }  i \le n  \\
 0 &\mbox{ else,}   
\end{array} \right .  \\ 
{\bar f}^n_b ({\x}, {\bar x}_b)& := & 
 \frac{1}{n}
 \sum_{j \le n}   \xi_{j} ({\bar x}_j, x_b) W_{j, b}  \ \mbox{ with } \  x_b     :=  f^b (  {\bar x}_b).
\nonumber
\end{eqnarray}} It is clear that the  fixed points   of the above operators  equal  the   aggregate vectors, $(\{{\bar X}_i^s \}_{i \le n} , {\bar X}_b )$.
Define the  'limit' operator $\bar {\bf f}^\infty (\x, {\bar x}_b)  =  ({\bar f}^\infty_b, {\bar f}^\infty_1, {\bar f}^\infty_2 \cdots )$:
\begin{eqnarray}
\label{Eqn_bar_fixedpoint_limtGen}
 {\bar f}^\infty_i ( \x, {\bar x}_b ) :=
\lim \sup_n {\bar f}^n_i (\x, {\bar x}_b) \mbox{ for all } i \in \{b, 1, 2, \cdots \}. 
\end{eqnarray}
The idea is to show that the fixed point of this operator equals that of   a  `limit' system and that the fixed points of the original system converge towards  these fixed points. 
Recall that the weights sum up to one, i.e., $\sum_{i} W_{j, i} = 1$ for any $i$.  
Thus we require the fixed point of the  operator:
$$
{\bar {\bf f}}^n   \mbox{ where }  {\bf f}^n  :  {[0, y]}   \times  {s}^\infty \to   {[0, y]}   \times  {s}^\infty,
$$
where  the $s^\infty$ is defined in footnote \ref{footnote_sinf}. 
Idea is to derive a kind of mean field analysis where the aggregates will be approximated by their expected values. 

When we consider constant sequence, i.e., if $\x = ({\bar x},  {\bar x},\cdots)$ the limit superiors in the definition of the  limit system ${\bar {\bf f}}^\infty$ are actually limits
by {\bf A}.2 and Law of large numbers (LLN) and equal (with $x_b$ as in (\ref{Eqn_bar_fixedpoint_randGen}))

\vspace{-8mm}
{\small \begin{eqnarray}
{\bar f}_i^\infty  (\x,  {\bar x}_b)  \hspace{-1.5mm}&  \hspace{-1.5mm}=  \hspace{-1.5mm}& \hspace{-1.5mm}     E_{G_i, \eta_i^{bs}} \left [ \xi_i ({\bar x}, x_b)   \right ]  (1-E[\eta_1^{sb}]) \mbox{ and }  \nonumber \\
{\bar f}_b^\infty  (\x, {\bar x}_b)  \hspace{-1.5mm}&  \hspace{-1.5mm}=  \hspace{-1.5mm}& \hspace{-1.5mm}     E_{G_i, \eta_i^{bs}}  [\xi_i ({\bar x}, x_b)  ] E[\eta_1^{sb}]  .  \label{Eqn_barf_limit}
\end{eqnarray}} 
In the above, $E_{X,Y}$ represents the  expectation with respect to  $X,Y$.
The random variables are IID and hence the first equation is the same function for all $i$.
By Theorem \ref{Thm_MainGen},  given below, 
one such constant sequence would be the almost sure limit of the solutions of the aggregate fixed point equations (\ref{Eqn_bar_fixedpoint_randGen}).   Thus one will have to solve a two-dimensional fixed point equation corresponding to the above function (\ref{Eqn_barf_limit}). And then random fixed points  (\ref{Eqn_FixedeqGen})-(\ref{Eqn_Fixedeq2Gen}) are asymptotically independent depending upon the other nodes only via  the aggregate fixed point, as given by the theorem below\ATR{.}{ (proof is in \cite{TR}).}

\begin{theorem}
\label{Thm_MainGen}
%Fix any realization of  $(G_c, G_b)$.  
Assume either  $0 < E[\eta_1^{sb}]   < 1$    or $\sigma < 1$ in {\bf A.3}.
The aggregates  of the random system,  which are   FPs of    (\ref{Eqn_xi})-(\ref{Eqn_bar_fixedpoint_randGen}), denoted by $({\X}^{*}, {\bar X}_b^*) (n) := ( \{{\bar X}_i^s\}_i, {\bar X}_b^* ) (n) $ converge   as $n \to \infty$ along a sub-sequence. That is there exist a $k_n\to \infty$ such that: 
\begin{equation}
\label{Eqn_main_aggregate_convergence}
 {\bar X}_i^s(k_n)  \to {\bar x}^{\infty*} \mbox{ for all } i \mbox{ and }  {\bar X}_b^*(k_n)  \to {\bar x}_b^{\infty*} \mbox{ almost surely (a.s.),}
 \end{equation}
 where $({\bar x}_b^{\infty*}, \x^{\infty*})$ with $\x^{\infty*} := ({\bar x}^{\infty*}, {\bar x}^{\infty *} , \cdots)$ is the FP of 
  the limit system given by  (\ref{Eqn_barf_limit}).
Further    (any sequence of) FPs of the  original system (\ref{Eqn_FixedeqGen})- (\ref{Eqn_Fixedeq2Gen}) converge almost surely (along the sub-sequence of  \eqref{Eqn_main_aggregate_convergence}, i.e., as $k_n\to  \infty$): %to the limit :
\begin{eqnarray}
X^b (k_n)  &\to&  X_b^{\infty* } :=  f^b (  {\bar x}_b^{\infty*}) \mbox{ as $n \to \infty$    and } \hspace{3mm}\label{Eqn_Act_Fixed} \\
X^s_i (k_n)  &\to& f^s (G_i,  {\bar x}^{\infty*}, \eta_i^{bs} X_b^{\infty*})  
.  \hspace{20mm}\mbox{ \eop } \nonumber
\end{eqnarray}
\end{theorem}
\ATR{
  {\bf Proof:}  This is a special case of   \cite[Theorem 1]{ANOR}  and the proof is available in \cite{ANOR}. \eop}{}\ATR{}{}
  
  Thus  the fixed points of the finite $n$ system converge to that of the limit system.  The fixed points are asymptotically independent  and depend upon the other nodes only via an almost sure constant ${\bar x}^{\infty*}_i$ which is common for all $i$. Another important point to observe here is that,  \underline{the aggregate fixed}  \underline{points need not be unique}, however any sequence of fixed points (one for each $n$) converges towards that of the limit system (when it has   unique fixed point).

 {\bf Remarks:}  Another interesting observation is that the result does not depend upon the precise probability $p_{ss}$ of the connection between small nodes.  It only depends upon the fact that every node can potentially influence every other node directly or indirectly (i.e., $p_{ss} > 0$).  
  It is straight forward to generalize  to the case   where\ATR{ $\{I_{i,j}\}$  are any IID random variables and  the weights  are formed in a similar way. Further}{} one may have finite number of groups of small nodes, nodes within a group are identical stochastically (identical $\{G_i\}$, $\{\eta_i^{sb}\}$ and $\{\eta_i^{bs}\}$), and any typical small node can be of group $i$ with probability $q_i$ independent of others.  With this one can study a wider variety of heterogeneous situations. For example one can consider a financial network with many big and small banks.  One can also  generalize to the case when we have more than one (but finite) distinct limits for  the aggregate fixed points.  The results are true even  when $\{X_i^s\}$ are  finite dimensional,   with dimension greater than one. One can then consider banks with different levels of connectivity.   %We now apply our results to a large Financial network.
  % For the rest of the paper we apply this result to an example financial network. 

\section{Financial Network}

  Consider a huge financial  network with $n$  small banks and one big bank. The assets (shares, bonds etc.) of the big bank are   large compared to any small bank and the small banks 
  are similar in nature. 
   At time  $T =0$    the banks invest in a project  by taking loans from one another or from outside the financial network. At time   $T=1$, the banks
anticipate some returns from their investments, which is used    
   to  clear their obligations (e.g.,   \cite{eisenberg2001systemic,acemoglu2015systemic}). But the investments are risky,  there are chances of economic shocks,  the returns might be lower than anticipated, because of which   some banks may not be able to (fully/partially) pay the liability.  We say these banks \underline{defaulted.} The defaulted banks increase the shocks to other connected banks, because of which we may have more defaults. And this can continue and the system can `collapse'.  Systemic risk precisely studies these aspects. 
\ForTR{The defaulted banks  break the bonds, those that they invested  at time $T=0$, and try to clear their obligations using the partial returns obtained after breaking.
At the end of second time period $T=2$ the survived banks obtain a return $A^s$  (respectively $nA^b$ for the big bank) while the defaulted ones obtain 
 obtain $\rho^s A^s / \rho^b n A^b$ (at time period $T=1$) where $\rho^s/\rho^b  < 1$ (usually much less than 1).}{} Our main aim is to study the influence of economic shocks on the stability of the network, wherein stability is understood in terms of the fraction of defaults, for a given  realization of the  shocks\ForTR{ or in terms of the expected surplus after the shocks etc.}{.} 
 
 \underline {Liabilities:} We model the financial network using a directed weighted graph, with the banks as the nodes and the weighted edges   represent the liability fractions and directions.  
The weight $W_{i,j} := l_{i,j} /  Y_i$  represents the fraction of the liability, where  $l_{i,j}$ is the amount the bank $i$ is liable to bank $j$ and   $Y_i := \sum_j l_{i,j} + l_{i,b}$ is the total liability of the small bank $i.$
The small banks are liable to  big bank with  proportionality factors $\{\eta^{sb}_i\}_i $ and,  a small bank is liable to another  small bank  with probability $p_{ss} >0$. Let $I_{i,j} $ be 1 if small bank $i$ is liable to   small bank $j$ and then the fractions of liability would be:
$$
\hspace{5mm}
W_{i, j} := \frac{I_{i,j} (1-\eta^{sb}_i) }{ \sum_{j'} I_{i,j'} }   \mbox{ and }     W_{i,b} = \eta^{sb}_i. \vspace{-0.5mm}
$$
The fraction of liability of  big   bank towards  small  bank $i$   equals $W_{b,i}=\eta^{bs}_i / {\bar \eta}$, \ForTR{and $W_{b,o} =\eta^o/ {\bar \eta}$ is the fraction that it is liable to sources outside the network,  with
 ${\bar \eta} := \sum_{j \le n } \eta^{bs}_j +  n \eta^o  $.}{ with ${\bar \eta} := \sum_{j \le n } \eta^{bs}_j .$}
 Let $nY^b $ represent the total liability of big bank.  

\underline{Shocks:}  Let  $K^s_{i}$ be the amount of money small bank $i$ is expecting as return (plus its liquid  cash) at time $T=1$, let 
 $Z^s_i$ represent the individual/independent  shock experienced by small bank $i$, and
let $Z_c$ represent the shock that is commonly received by all the small banks.  The big bank receives a  shock of magnitude $n \delta Z_c$ along with its independent shock 
$n Z_b$.  After the shocks the small bank $i$ receives  
    $(K^s_{i}- Z_c- Z^s_{i})^+$ at $T=1$ while the big bank receives $n(K^b - \delta Z_c- Z^b )^+$, where   $nK^b$ is the shock free return anticipated at time $T=1$. 
\begin{figure*} 
 \vspace{-5mm}
  %      \begin{minipage}{10.6cm}
       \centering
       \hspace{-5mm}
              \begin{subfigure} 
               \centering
               \includegraphics[width=5.cm,height=6cm]{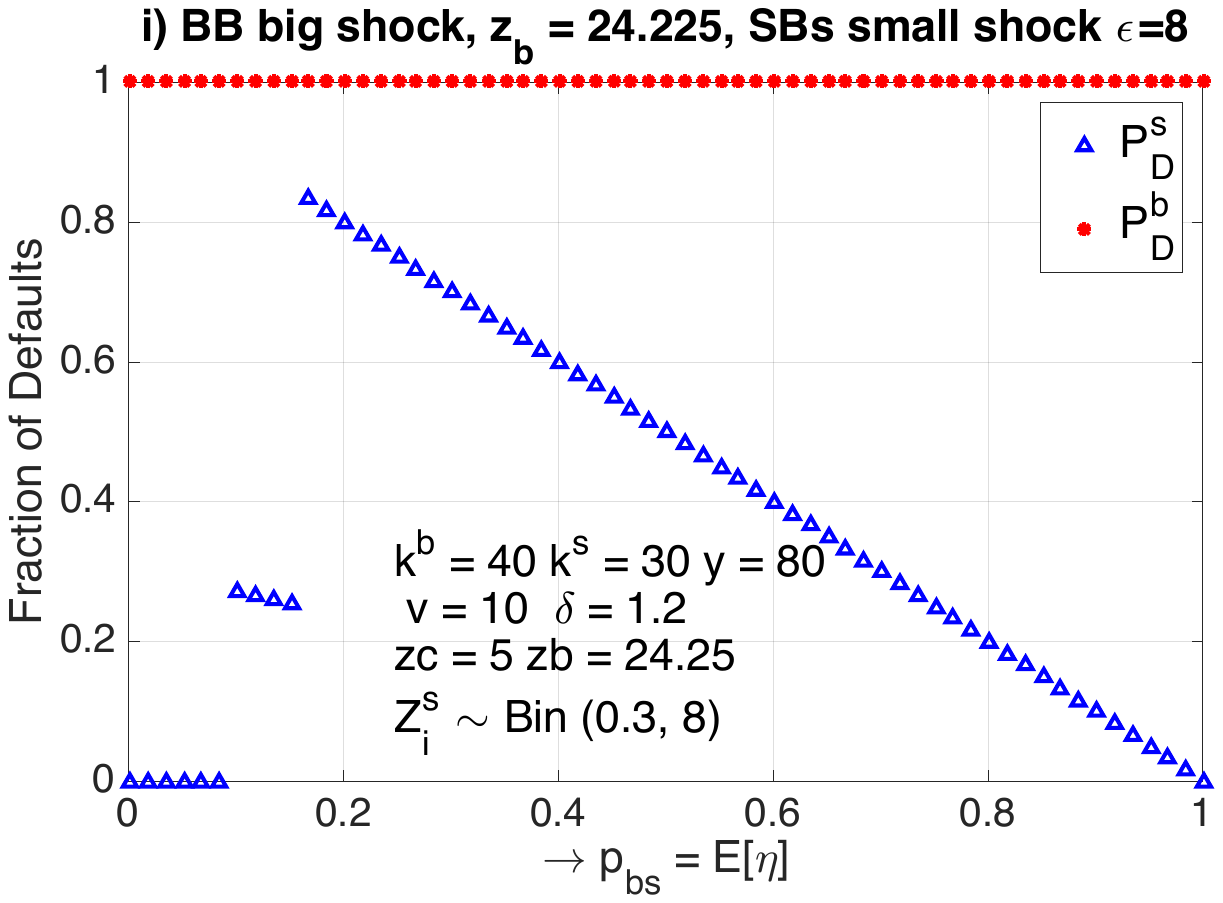}
            %
            %  \caption{Expected fraction of defaults}
       \end{subfigure} 
\hspace{.1mm}
       \begin{subfigure} 
               \centering
               \includegraphics[width=4.9cm,height=6cm]{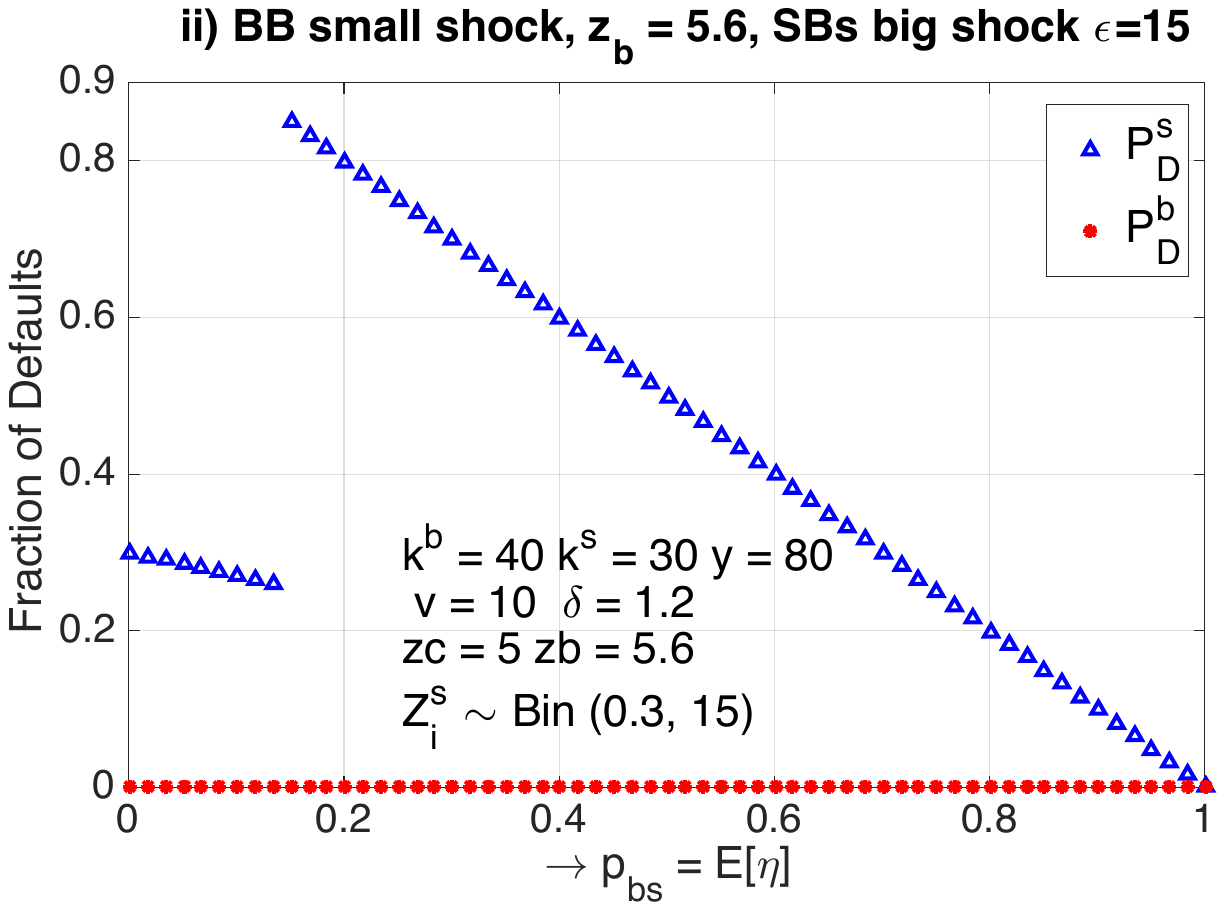}
            %  \caption{Expected fraction of defaults}
       \end{subfigure}   \hspace{.1mm}
              \begin{subfigure}
               \centering
               \includegraphics[width=5.cm,height=5.7cm]{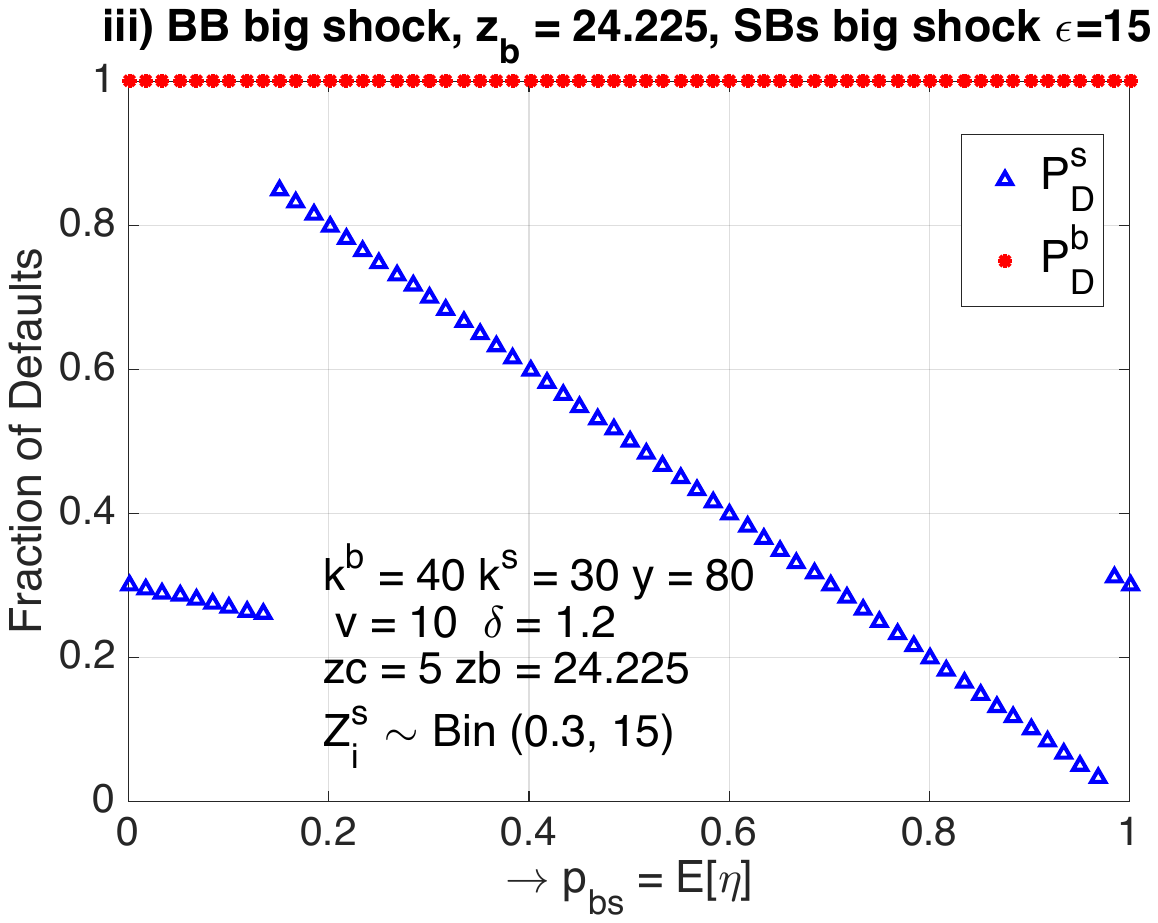}
               %\caption{Expected Surplus till time $T=1$}
             %  \caption{Big bank fails}
       \end{subfigure} 
       \vspace{-8mm}
       \caption{ First sub-figure: only BB defaults;  Second sub-figure: only SBs default;  Third sub-figure: both default. 
        \label{Fig_All_retuns} }
    \end{figure*}

\ForTR{

\begin{figure*} 
 \vspace{-4mm}
  %      \begin{minipage}{10.6cm}
       \centering
       \hspace{-4mm}
       \begin{subfigure}
               \centering
               \includegraphics[width=5.2cm,height=8cm]{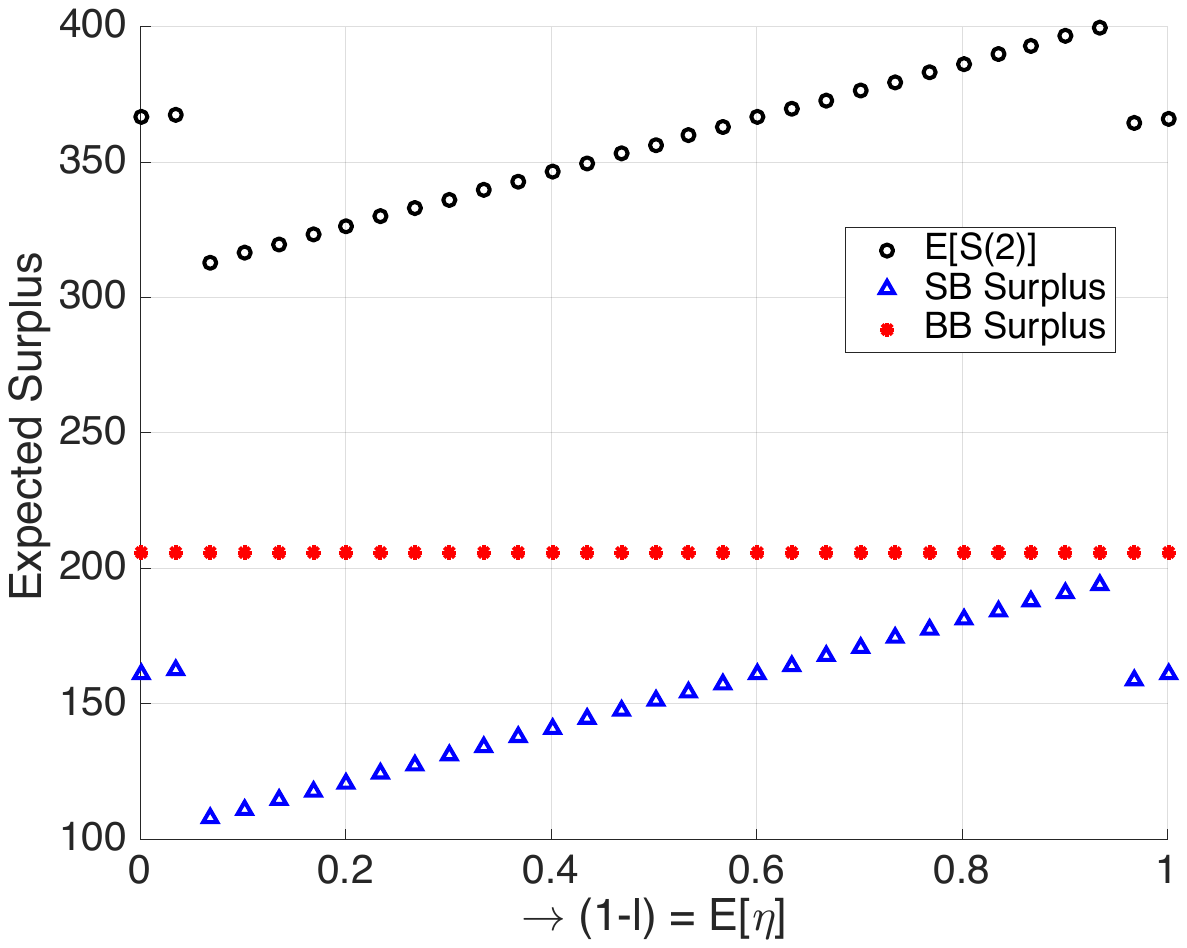}
               %\caption{Expected Surplus till time $T=1$}
       \end{subfigure} \hspace{-1mm}
       \begin{subfigure} 
               \centering
               \includegraphics[width=5.2cm,height=8cm]{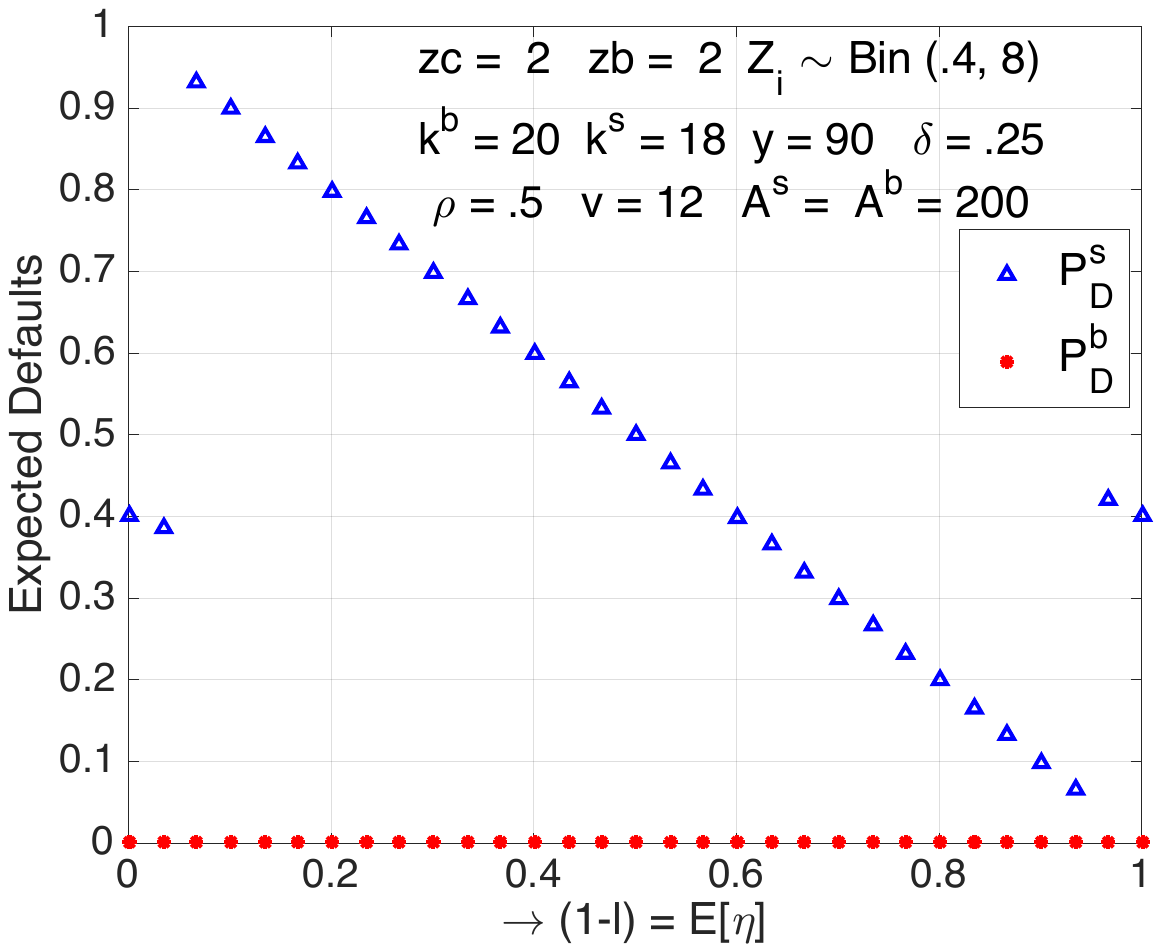}
            %  \caption{Expected fraction of defaults}
       \end{subfigure} 
       \vspace{-12mm}
       \caption{ Small bank-major shocks ($Z_i \sim Bin (0.4, 8)$), big bank-minor shocks ($Z_c = 2, Z_b= 2$)  connection with big bank (best at 0.9778) improves the surplus.
        \label{Fig_SmallSB_retuns} }
\end{figure*} 
 }{}

\underline{Clearing Vectors:}
 Let $X_i^s$ represent the maximum possible part of the total liability, eventually cleared by small bank $i$, and let $X_b$ (per small bank) represent the same for big bank. 
 The vector $(X^b,  \{X_i^s \}_s)$ is referred to  as clearing vector (e.g., \cite{eisenberg2001systemic,acemoglu2015systemic}) and we make the following commonly made assumptions (as in \cite{eisenberg2001systemic,acemoglu2015systemic} etc.) for computing the same.

  When a bank (say bank $i$)  defaults, it may not be able to clear its liabilities completely. However  it repays   the maximum  possible, and the amount cleared to another  bank  (say bank $j$) is proportional to the fraction $W_{i,j}$. 
  Thus small bank $i$ receives ${\bar X}_i^s := \sum_{j} W_{j,i} X_j^s$  at time period $T=1$ when the other small banks try to clear their liabilities. In a similar way, it receives 
  $X^b W_{b, i}$ from big bank.  
  It also   receives  $(K^s_{i}- Z_c- Z^s_{i})^+$ at time period $T=1$ (after shock) from outside investments.   The liabilities are paid, only after paying the operational costs/taxes $v^s$ ($nv^b$ for big bank).  Thus the bank $i$  at maximum can clear $( (K^s_{i}- Z_c- Z^s_{i})^+ + {\bar X}_i^s - v^s)^+$\ForTR{ and if this amount is less than $Y_i$ it breaks its bonds which are supposed to mature at time period $T=2$.}{.} 
The amount cleared by big bank  is also computed in a similar manner. 
  Thus, in all, the total amount cleared by small bank $i$ and big bank  respectively equals,
 
 \vspace{-6mm}
{\small
\begin{eqnarray*}
X_i^s  \hspace{-1.5mm}&  \hspace{-1.5mm}=  \hspace{-1.5mm}& \hspace{-1.5mm} \min  \Big \{ \Phi_i^s ({\bar X}^s_i, X^b) ,   
    Y_i \Big \}   ,  \mbox{ and }\\
  X^b \hspace{-1.5mm}&  \hspace{-1.5mm}=  \hspace{-1.5mm}& \hspace{-1.5mm} \frac{1}{{\bar \eta}} \min \left \{    n \Phi^b ({\bar X}^b)  ,  \    nY^b \right \}   \mbox { with  aggregates }   
    \\
{\bar X}^b  \hspace{-1.5mm}&  \hspace{-1.5mm}=  \hspace{-1.5mm}& \hspace{-1.5mm}  \frac{  \sum_{j \le n} X_j^s W_{j,b} }{n}  , \ {\bar X}_i^s := \sum_{j} W_{j,i} X_j^s ,  \mbox{ where} \\
  \Phi_i^s ({\bar X}^s_i, X^b) \hspace{-1.5mm}&  \hspace{-1.5mm}:=  \hspace{-1.5mm}& \hspace{-1.5mm} \left (   (K^s_{i}-Z_c - Z^s_{i})^++ {\bar X}_i^s   +\ForTR{ \rho^s A^s+ }{} \eta_{i}^{bs} X^b - v^s\right )^+,\\
   \Phi^b ( {\bar  X}_b) \hspace{-1.5mm}&  \hspace{-1.5mm}:=  \hspace{-1.5mm}& \hspace{-1.5mm} 
  \left (   (K^b-Z_c\delta- Z^b)^+ \ForTR{+  \rho^b A^b}{} - v^b    +  {\bar X}^b \right )^+ .
  \end{eqnarray*}}
%  For   convenience, we redefine $X^b/{\bar \eta}$ as $X^b$ and    also redefine the other terms as below:
%  \begin{eqnarray*}
%  \Phi_i^s \hspace{-1.5mm}&  \hspace{-1.5mm}=  \hspace{-1.5mm}& \hspace{-1.5mm}   (K^s_{i}-Z_c - Z^s_{i})^++ {\bar X}_i^s   + \ForTR{\rho^s A^s+}{}   \eta_{i}^{bs} X^b - v^s \\
%  %
%    X^b  \hspace{-1.5mm}&  \hspace{-1.5mm}:=  \hspace{-1.5mm}& \hspace{-1.5mm}  \frac{n}{{\bar \eta}}    \min \left \{      \left ( \Phi^b \right )^+ ,  \   Y^b \right \} .
%\end{eqnarray*}    
By Law of large numbers (LLN),  
${\bar \eta} / n \to E[\eta^{bs}] \ForTR{+ \eta^o}{}$   a.s.  
The rest of the system is exactly like the general system   (\ref{Eqn_FixedeqGen})-(\ref{Eqn_Fixedeq2Gen}), 
for any given realization of $(Z_c, Z_b, K^b, Y^b).$  Assumptions {\bf A.1} and {\bf A.3}  are clearly satisfied and Theorem \ref{Thm_MainGen} is applicable if we assume {\bf A.2}.  
By  Theorem \ref{Thm_MainGen} and equation (\ref{Eqn_barf_limit}), for any given realization $(Z_c, Z_b, K^b, Y^b) = (z_c, z_b, k^b, y^b)$, 
the aggregate vectors are approximately (accurate for large $n$) the solutions of  
  the following   fixed point equations: 
 \begin{eqnarray}
 {\bar x}^{\infty*} \hspace{-1.5mm}&  \hspace{-1.5mm}=  \hspace{-1.5mm}& \hspace{-1.5mm}   
  E_{Z_i, Y_i, K_i, \eta_i^{sb}} \hspace{-.5mm} \Big [ \hspace{-1.mm}  \min \Big \{  \Phi^{s}_i ({\bar x}^{\infty*} \hspace{-.5mm}, { x}^{\infty *}_b),  Y_i   \Big  \} \Big ]  \hspace{-.5mm}(1 - E[\eta_1^{sb}]) , \nonumber \\
  %
%\Phi^{s\infty}_i   \hspace{-1.5mm}&  \hspace{-1.5mm} :=  \hspace{-1.5mm}& \hspace{-1.5mm}
%  (K^s_{i}-z_c- Z^s_{i})^++    {\bar x}^{\infty*}   \ForTR{+  \rho^s A^s}{}+   { x}^{\infty*}_b  \eta^{bs}_i - v^s,   \nonumber \\
  %
   {\bar x}^{\infty *}_b  \hspace{-1.5mm}&  \hspace{-1.5mm}=  \hspace{-1.5mm}& \hspace{-1.5mm}  
      {\bar x}^{\infty*}  \frac{E[\eta_1^{sb}]}{ 1 - E[\eta_1^{sb}]}, \ \  %\nonumber  \\
    { x}^{\infty *}_b  = %\hspace{-1.5mm}&  \hspace{-1.5mm}=  \hspace{-1.5mm}& \hspace{-1.5mm} 
       \frac{ \min \left \{   \Phi^b ({\bar x}^{\infty *}_b), 
  %(k^b - z_c \delta - z_b)^+ +    {\bar x}_b^{\infty*} \ForTR{+ \rho^b A^b}  - v^b   ,
    \ y^b \right \}  }{E[\eta_1^{bs}] \ForTR{+ \eta^o}{}} . \label{Eqn_Finanace_FP}
 \end{eqnarray} 
 Once these fixed point equations are solved, the clearing vectors are approximately given by (\ref{Eqn_Act_Fixed}) of Theorem \ref{Thm_MainGen}.  These are \underline{asymptotically independent}, as now the aggregates ${\bar X}^s_i$ are almost sure constants and are common for all $i$.  We now compute relevant asymptotic performance measures. %Thus one can easily compute interesting performance measures. % as discussed below. 
\vspace{-4mm}

\ForTR{

\subsection{Performance measure to study Systemic risk}
Once these fixed point are available for each pair of shock, initial value, total liability realizations $(k^b, y^b, z_c, z_b)$,  one can obtain various performance measures as below:

1) \underline{Expected surplus till time $T=1$ per small bank:} This is the total income of the network (big bank as well as small banks) after paying away the liabilities and taxes $(v)$ divided by number of small banks.  We are currently using the following expression which has to be proved as in \cite{eisenberg2001systemic,acemoglu2015systemic}
\begin{eqnarray}
\label{Eqn_ES1}
E [S(1)] \lesseq   E \left [  \left ( \Psi^s  \right )^+  +  \left (\Psi^s +   \rho^s A^s  \right )^+ ;  \Psi^s < 0  \right ]  \nonumber \\ 
   &&  \hspace{1mm} +  (\Psi^b )^+ +     \left (\Psi^b + \rho^b A^b  \right )^+  1_{ \{  \Psi^b < 0   \} },
  \mbox{ with} \nonumber \\
    \Psi_i^s  \lesseq
  (K_i^s -Z_c- Z^s_{i})^++    {\bar x}^{\infty*}  +     { x}^{\infty*}_b    \eta^{bs}_i -v^s  -  Y_i^s \nonumber \\
    \Psi^b \lesseq   (K^b - Z_c  \delta - Z_b)^+ +   {\bar x}^{\infty*}   \frac{E[\eta_1^{sb}]}{ 1 - E[\eta_1^{sb}]}  \nonumber \\ && \hspace{10mm}  -v^b   -  Y^b .
 \end{eqnarray}
 The expectations are with respect to $(Z_i, \eta_{i}^{bs}, Y_i, K_i^s)$ and are   conditioned over $(K^b, Z_c, Z_b, Y^b)$.

 2) \underline{Expected number of defaults:}}
 {\underline{Fraction of defaults:}}  The  fraction of small banks that defaulted  (by bounded convergence theorem) and the indicator that the big bank defaults
asymptotically equal:
 
 \vspace{-8mm}
 {\small
\begin{eqnarray}
\label{Eqn_PDs}
\hspace{-2mm}
P_D^s \hspace{-1.mm} := \hspace{-1.mm} P  \left ( \PhiPsi_i^s ({\bar x}^{\infty*}  , { x}^{\infty *}_b)  <  Y_i   \right ),      \  
P_D^b  \hspace{-1.mm} := \hspace{-1.mm}   1_{ \{ \PhiPsi^b (  {\bar x}^{\infty *}_b)  < y^b   \} } \mbox{a.s.}   \hspace{-1mm}
 \end{eqnarray}} 
 \ForTR{
 3) \underline{Expected surplus till time $T=2$ per small bank:}  The banks invest and the liability structure is defined at time period $T=0$. The first installment is returned  to the banks at $T=1$ (as in \cite{eisenberg2001systemic,acemoglu2015systemic}) and we are studying the situation when there are shocks to these returns.  The surplus per small bank till time period $T=1$ is given by $E[S(1)]$. 
 The banks that defaulted at time $T=1$, break their bonds/investments. If this did not happen a small bank receives amount $A^s$ on maturity at $T=2$, while the big bank   receives $nA^b$ amount (if it did not default). Thus the surplus at $T=2$ equals:
 \begin{eqnarray}
 \label{Eqn_ES2}
E[S(2)] =   E[S(1)]  + (1-P_D^s ) A^s + (1-P_D^b) A^b . \hspace{2mm}
 \end{eqnarray}

}{}

 \ForTR{ 
\begin{figure*} 
 \vspace{-22mm}
       \centering
       \hspace{-9mm}
       \begin{subfigure}
               \centering
               \includegraphics[width=5.9cm,height=8cm]{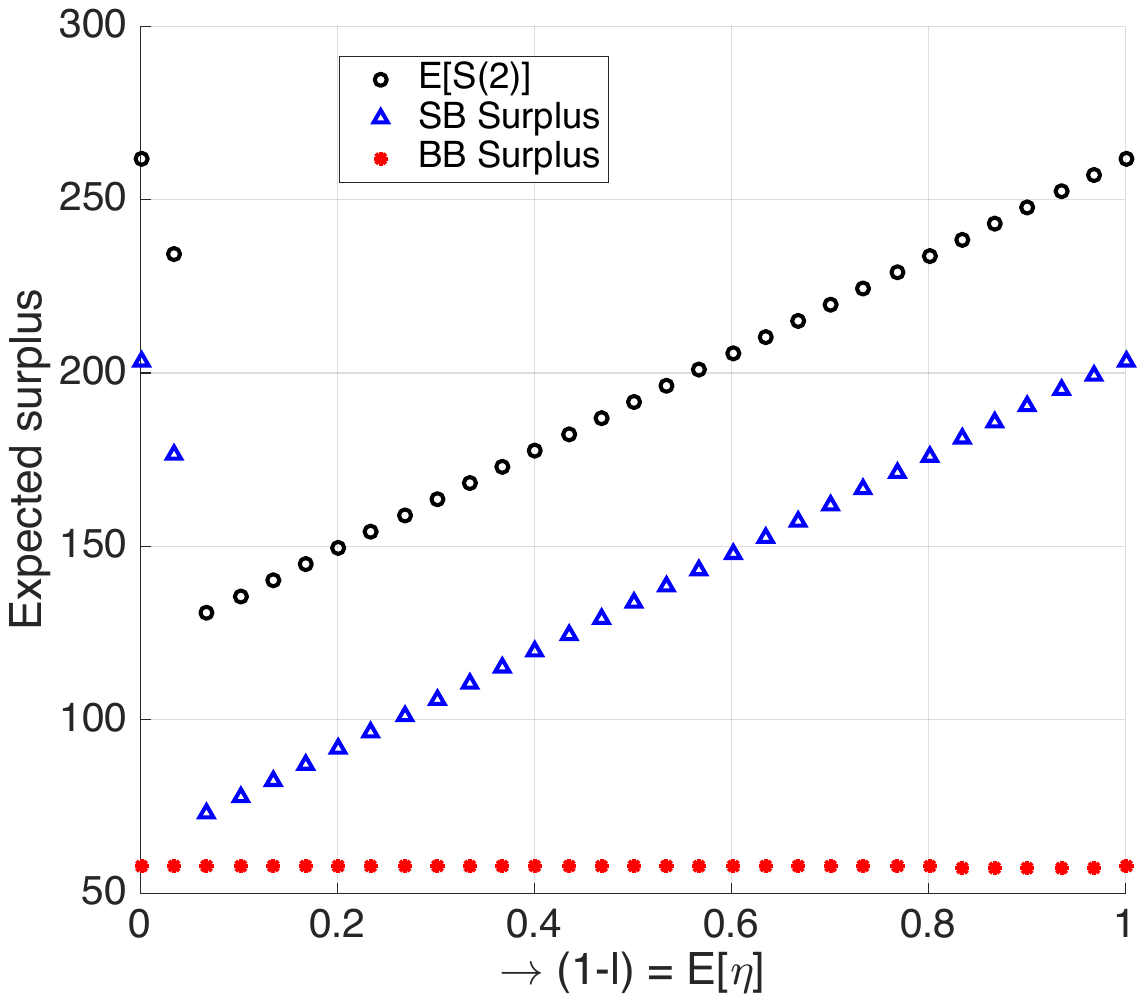}
               %\caption{Expected Surplus till time $T=1$}
       \end{subfigure} \hspace{3mm}
       \begin{subfigure} 
               \centering
               \includegraphics[width=5.9cm,height=8cm]{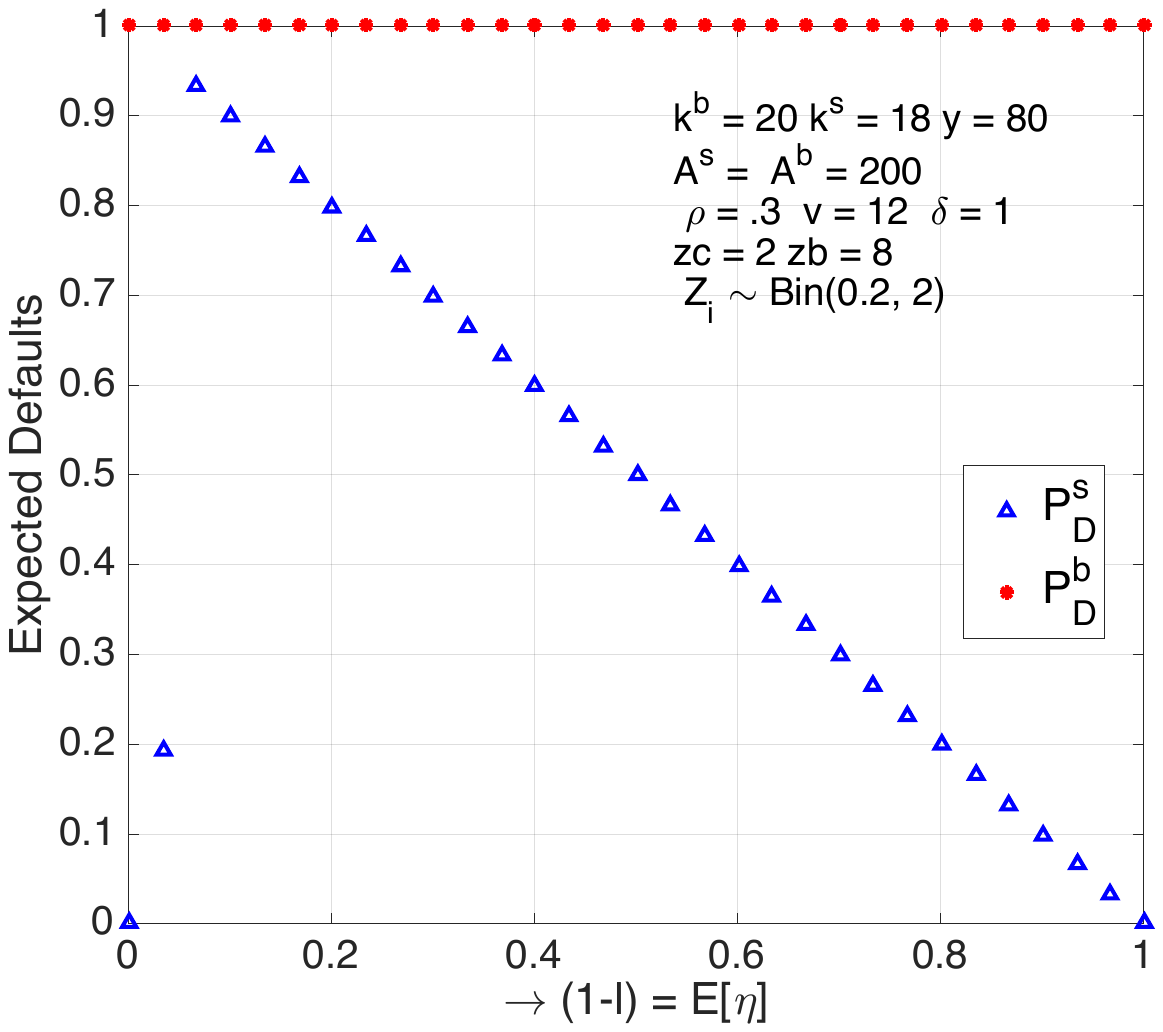}
            %  \caption{Expected fraction of defaults}
       \end{subfigure} 
       \vspace{-20mm}
       \caption{ When  big bank has major shocks ($Z_c = 2, Z_b= 8$),  risk of small bank is minor  ($Z_i \sim Bin (0.2, 2)$). Connection with big bank does not help, 
       surplus at $p_{bs} = 0$ equals  that at $p_{bs}=1$.  
        \label{Fig_Big_BB} }
   %     
 %       \end{minipage}  
        %
        %
 
        %
        %
%        \begin{minipage}{10.6cm}
%
       \centering
       \vspace{-9mm}
       \begin{subfigure}
               \centering
               \includegraphics[width=5.2cm,height=8cm]{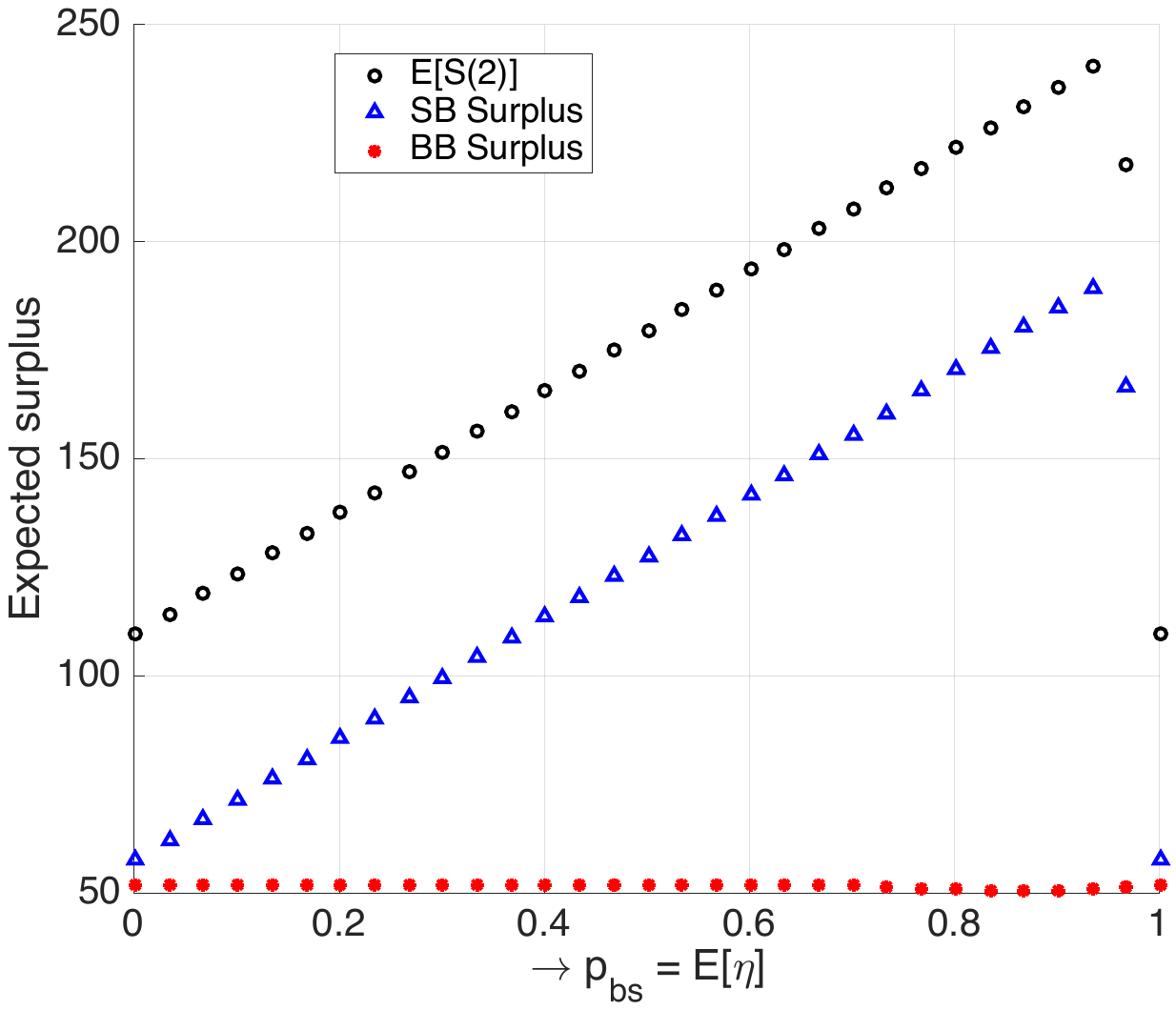}
               %\caption{Expected Surplus till time $T=1$}
       \end{subfigure} \hspace{-1mm}
       \begin{subfigure} 
               \centering
               \includegraphics[width=5.2cm,height=8cm]{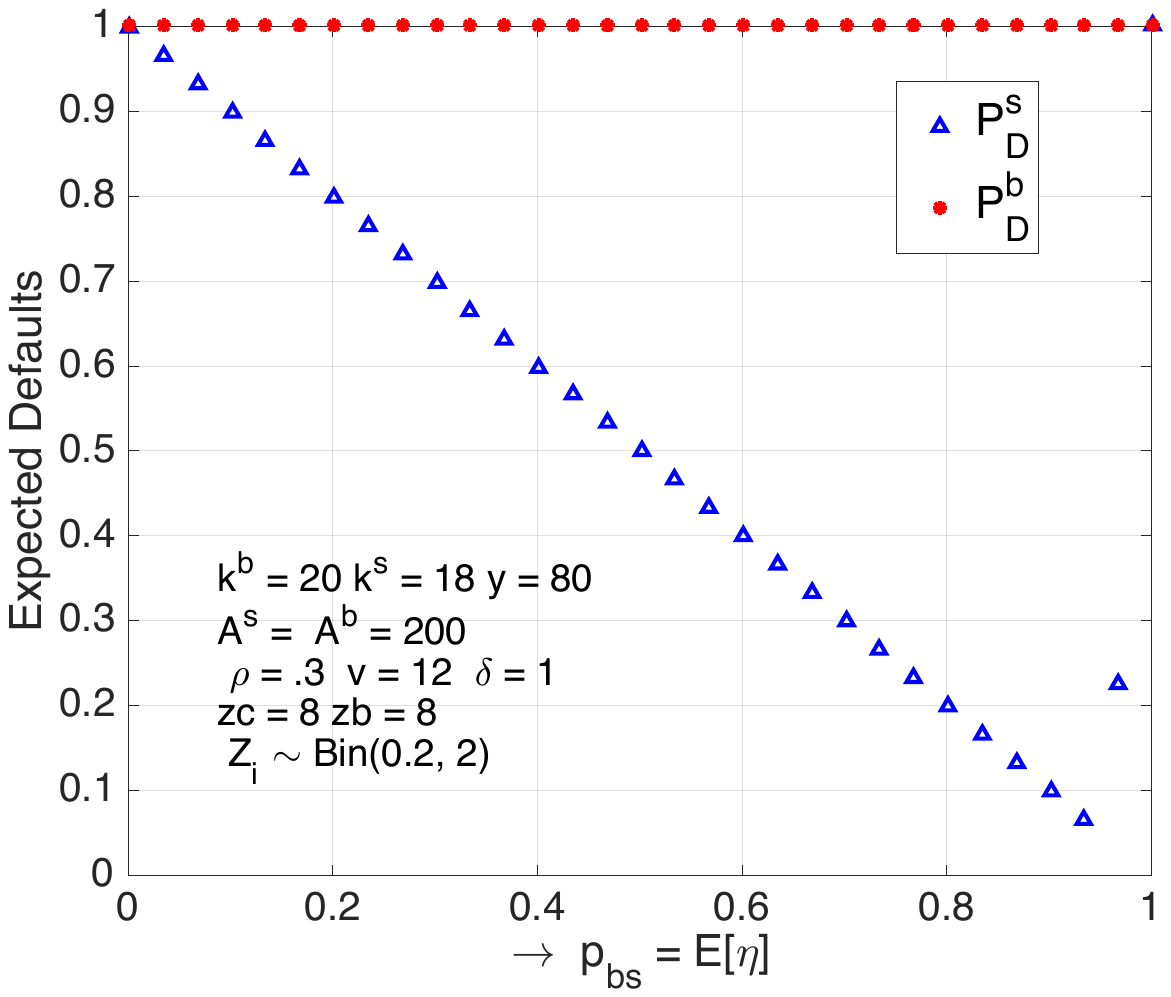}
            %  \caption{Expected fraction of defaults}
       \end{subfigure} 
       \vspace{-5mm}
       \caption{ With big common shock, setting as in Figure \ref{Fig_Big_BB} but for $z_c = 8$.  Now the connection with big bank  improves the performance as in Figure \ref{Fig_SmallSB_retuns}. The improvement is more pronounced, but big-bank performance remains the same. 
        \label{Fig_ZcZbBig} }
        
  %      \end{minipage}
    
         \vspace{-5mm}
\end{figure*}}{}  

\underline{Regular networks:}  In these   networks,  the total claim of any bank equals  its total liability, i.e.,  $ \sum_{j} l_{ji}+l_{bi} = \sum_{j} l_{ij} +l_{ib}  $ for all $i$.
We consider these networks for further study, as they ensure the initial  wealth\footnote{This is proportional to the amount anticipated without shocks at time $T=1$ plus the amount anticipated by the returns from other banks minus the amount it has to pay to other banks, all at time $T=1$.  
For small banks and big bank (per small bank) it is proportional respectively to: 
 
\vspace{-6mm}
{\small $$
K^s_i + \sum_{j} l_{j, i} + l_{b,i}  - \sum_{j} l_{i,j} - l_{i, b}  \mbox{ and }  K^b +  \sum_{j} l_{j, b}    - \sum_{j} l_{b,j}.
$$}} of the network remains the same
once the characteristics of $K^b, K^s, Z_i, Z_c, Z_b$ remain  the same (stochastically).  
 This allows fair comparison of  stability of the  network for different values of the parameters\ForTR{, importantly the big bank-small bank connection parameter $E[\eta^{sb}]$. }{.}
 In our network  the liabilities are  random and equal:

\begin{eqnarray*}
l_{ji}  = Y_j \frac{I_{ji} (1-\eta^{sb}_j) }{\sum_{i'} I_{j,i'}} \mbox{,  } \  l_{jb}  = \eta^{sb}_j Y_j   \mbox{ and }  
l_{bj}=  \frac{1}{{\bar \eta}}  \eta^{bs}_j nY^b .%  \mbox{ and thus } \\
%\sum_{i} l_{ji} &=& \sum_i   Y_j \frac{I_{ji} (1-\eta^{sb}_j) }{ |A_j|} = Y_j  (1-\eta^{sb}_j) . 
\end{eqnarray*}
Thus clearly  $\sum_{i} l_{ji} +  l_{jb} = Y_j$, and  as $n \to \infty$ by LLN and {\bf A.2}.

\begin{eqnarray*}
\sum_{i} l_{ij}   \hspace{-1.mm}&  \hspace{-1.mm}=  \hspace{-1.mm}& \hspace{-1.mm} \sum_i   Y_i \frac{I_{i,j} (1-\eta^{sb}_i) }{ \sum_{j}I_{i,j} } \to E[Y_i]  E[ 1-\eta^{sb}_i ] \mbox{,  thus} \\
\sum_{i} l_{ij} + l_{bj}  \hspace{-1.mm}&  \hspace{-1.mm}  \to  \hspace{-1.mm}& \hspace{-1.mm}   E[Y_i] (1 - E[\eta^{sb}_i ] )+  \frac{1}{E[\eta^{bs}] \ForTR{+ \eta^o}{}}  \eta^{bs}_j Y^b  .
\end{eqnarray*} 
Since the liabilities are random, we require equality in  stochastic sense or atleast in expected sense ($\stackrel{m}{=}$), i.e., we need:
{\small $$
 E[Y_i] (1 - E[\eta^{sb}_i ] )+  \frac{1}{E[\eta^{bs}]  \ForTR{+ \eta^o}{}}  \eta^{bs}_i Y^b  \stackrel{m}{=}  Y_i
$$}
We also need that  the liabilities and claims  of the big bank match, i.e.,
$
\sum_j l_{bj}  = \sum_j l_{jb} .
$ That is we need (at limit),
{\small $$
\lim_{n \to \infty} \frac{\sum_j l_{bj}}{n}  = \lim_{n \to \infty}  \frac{ \sum_j \eta^{bs}_j  }{nE[\eta^{bs}] \ForTR{+\eta^o}{}  } Y^b   
\  \stackrel {m}{=} \  
%\frac {\sum_j l_{jb} }{n}  = 
\lim_{n \to \infty} \frac{ \sum_j  \eta^{sb}_j Y_j }{n} .
$$}

\vspace{-7mm}
\subsection{Example Case studies}

\ForTR{
We consider two scenarios of regular networks and compute the performance measures (\ref{Eqn_ES1})-(\ref{Eqn_ES2}).

\subsection{No external links, $\eta^o = 0$}}{
We consider an example scenario of a regular network  and compute the performance measures.
}To keep things simple  yet sufficiently interesting, we consider a deterministic $  Y_i^s \equiv y, Y^b \equiv y^b $ and $K_i^s \equiv k^s$.  We then consider a given scenario 
$(z_c, z_b, k^b)$ as discussed before. This is the case with identical small banks in terms of initial wealth, investments at time $T=0$ and when they receive common shock $z_c$ 
as well individual independent shocks $\{Z_i^s\}$. The total shock of the big bank equals $\delta z_c +z_b$. To have regular networks  we set: 
$$
\hspace{22mm}\eta^{bs}_i  \stackrel{d}{=}  \eta^{sb}_i  \mbox{ and }  y^b   =   y  E[\eta^{bs}]  \mbox{ for all } i,
$$and so $E[\eta_{i}^{sb}] =E[\eta_{i}^{bs}]= p_{bs}$. 
\ignore{We also need that  the liabilities and claims  of the big bank match, i.e.,
$$
\sum_j l_{bj}  = \sum_j l_{jb} 
$$But this is readily satisfied since, 
$$
\sum_j l_{bj}  =  \frac{1}{E[\eta^{bs}]  }  \sum_j \eta^{bs}_j Y^b    \mbox{ and } \sum_j l_{jb}  = \sum_j  \eta^{sb}_j Y_j .
$$

The fixed point equations for this special case are the following. For any given realization of $(z_c, z_b)$ the common shock and big bank shock, the liability vectors in the limit are given by solving:
\begin{eqnarray*}
 {\bar x}^{\infty*} (z_c, z_b) &=&  \boxed{ E_{Z_i^s, \eta^{bs}_i} \left [  \min \left \{   
  (k^s -z_c- Z^s_{i})^++    {\bar x}^{\infty*} l  +     {\bar x}^{\infty*}_b  \eta^{bs}_i -v,  \ y \right \} \right ]} \\
    {\bar x}^{\infty *}_b (z_c, z_b) &=& \frac{1}{E[\eta^{bs}_1]} \boxed{   \min \left \{   
  (k^b - z_c  \delta - z_b)^+ +   {\bar x}^{\infty*}(1- l) -v    ,  \ y^b \right \} }
 \end{eqnarray*} with $y^b = y(1-l)$ and $E[\eta^{bs}_1] = (1-l)$ for regular networks.
}We immediately have the following   for the limit system  when $\eta_{i}^{bs}$ are indicators. %,  with {\small $E[\eta^{bs}] = p_{bs}$}.  %Note that for regular networks $E[\eta_{i}^{sb}] = p_{bs}$. 
Let $\Kmin$, $\Kmax$ respectively represent the worst and best returns ($(k^s_i-z_c-Z_i^s)^+$) of a small bank, given $z_c$. Then 
\ATR{}{\vspace{-5mm}}
\begin{lemma}
\label{Lemma_First} % Let $y > \Kmax-\Kmin.$
(i) If $y p_{bs} \le  \big (\Kmin    - v^s   \big ) $  then none of the small banks default, i.e, {\small $P_D^s = 0$}.  \\
(ii) If  {\small $yp_{bs} > \big  ( \Kmax +    x^{\infty *}_b - v^s \big )  $}, then  all   small banks default, i.e., {\small $P_D^s = 1$}. Thus if 
{\small $y p_{bs} > \big  ( \Kmax  +  y - v^s \big )  $}, then {\small $P_D^s = 1$}.
\ForTR{\\ (iii) If  {\small $0 < yp_{bs} < \big  ( \Kmax +    x^{\infty *}_b - v^s \big ) $}, then atleast some small banks do not default,  as $P_D^s \le 1- (1-w)pbs < 1$.}{}
\end{lemma}
{\bf Proof:}  The  small banks never default if for all  scenarios (realizations of $Z_i, \eta_i^{bs}$):

\vspace{-6mm}
  {\small $$ (k^s - z_c - Z_i)^+ + y (1-p_{bs}) - v^s + \eta_i^{bs} x_b \ge y. $$}
  The worst scenario is with  worst shock $\Kmin$ and with $\eta_i^{bs}=0$ and hence $P_D^s = 0$ when  $y p_{bs} \le \Kmin - v^s $  proving (i). In a similar way consider the best scenario   to obtain part (ii).  \eop
  
Thus we identified the conditions for zero and all defaults.  
As long as  $p_{bs} < (\Kmin-v^s)/y$, none of the small banks default.  But (for example) when $p_{bs}$ increases beyond $(\Kmin-v^s)/y$,   there can be a   `phase transition'  in 
the fraction of defaults,   $P_D^s$. At this point it probably would jump from $0$ to some non-zero value. One need more analysis to understand this possible `phase transition'.
We derive more such details for the special case with binary shocks. 

Consider that  {\small $Z_i  \sim  Bin (w, \epsilon)$}, i.e.,   binary $(0, \epsilon)$ shocks with {\small $P(Z_i = \epsilon) = w$}. 
    Then {\small $\Kmin =  (k^s - z_c -\epsilon)^+$} and {\small $\Kmax = (k^s -z_c)^+$}. 
   \ForTR{
 \subsubsection*{With two time periods}    
 We begin with analysis with two time period, $T = 0, 1$. Thus $\rho^s\ \rho^b$ and $A^s \ A^b$ are not applicable.  We compute only the expected fraction of defaults.
   }{}% 
   We have the closed form expressions for the clearing vectors  as well as the asymptotic fraction of defaults, for the sub-case when the big bank does not default. These expressions approximately equal the corresponding quantities  for system with large number of small banks. 
   
   \ForTR{
     \subsection*{When the big bank does not default}}{}
\begin{lemma}
\label{Lemma_BinaryNoRho}
Let  $y > \epsilon.$   With binary shocks,  the a.s.  limit fraction of defaults equal (with {\small ${\bar K}^s_Z := E[(k^s-z_c-Z_i)^+]$}):

\vspace{-8mm}
{\small    \begin{eqnarray*}
\hspace{-3mm}
P_D^s (p_{bs})\hspace{-1.2mm} &\hspace{-1.2mm} =\hspace{-1.2mm} &\hspace{-1.2mm}  \left [ \hspace{-2mm} \begin{array}{lllll} P_{D1} \\ 
P_{D2} \\
P_{D3} \\
P_{D4} \\
P_{D5} \end{array} \hspace{-2mm}\right ]  \hspace{-1.2mm}  = \hspace{-1.2mm}  \left \{ \hspace{-1.2mm}
\begin{array}{lllll}
0           & \mbox{if } \  b_0< yp_{bs}  \le   b_1  \\
w (1-p_{bs}) \hspace{-1.2mm}  & \mbox{if } \  b_1 <yp_{bs} \le b_2\\
1-p_{bs}            & \mbox{if } \  b_2 < yp_{bs} \le   b_3   \\
1- p_{bs}(1-w)  \hspace{-1.2mm}      & \mbox{if } \  b_3  < yp_{bs} \le   b_4\\
1                          & \mbox{else}  ,
\end{array} \right . \mbox{{\normalsize with,} } \\
&& \hspace{-15mm}b_i\   := \ c_i  (1-p_{bs}) P_{Di} + d_i  (1 - (1-p_{bs}) P_{Di}),  i > 0, \ b_0 = 0,
\end{eqnarray*}} 
%where the bounds are defined as below:
   %
%    \vspace{-3mm}
%{\small        \begin{eqnarray*}
%b_i &=& c_i  (1-p_{bs}) P_D^s + d_i  (1 - (1-p_{bs}) P_D^s)  
%%    b_1 &=& \Kmin  - v^s ,  \\ 
%%     b_2 &=&  c(P_D^s)  (1-p_{bs}) P_D^s + (\Kmax - v^s) (1 - (1-p_{bs}) P_D^s) , \\  
%%    b_3 &=& c(P_D^s)  (1-p_{bs}) P_D^s + (\Kmin - v^s + y) (1 - (1-p_{bs}) P_D^s),\\  
%%    b_4 &=& c(P_D^s)  (1-p_{bs}) P_D^s + (\Kmax - v^s + y) (1 - (1-p_{bs}) P_D^s) ,
%    \end{eqnarray*}}% using the constants  (with ${\bar K}^s_Z := E[(K^s_i-Z_c-Z_i)^+]$)
%    
    \vspace{-7mm}
{\small    \begin{eqnarray*}
\left [ \hspace{-2mm} \begin{array}{lllll} c_1 \\ 
c_2 \\
c_3 \\
c_4 \end{array} \hspace{-2mm}\right ]   \hspace{-1mm} =\hspace{-1mm} \left [\hspace{-2mm}
\begin{array}{lllll}
0            \\    
 \Kmin-v^s    \\
 {\bar K}^s_Z  - v^s \\
  \frac{ {\bar K}^s_Z   (1-p_{bs}) + (\Kmin +y) w p_{bs} }{1- p_{bs}(1-w) } - v^s  \hspace{-1mm}
\end{array}  \hspace{-1mm} \right ], 
\left [ \hspace{-1mm} \begin{array}{lllll} d_1 \\ 
d_2 \\
d_3 \\
d_4 \end{array}  \hspace{-1mm} \right ]   \hspace{-1mm}  = \hspace{-1mm} \left [\hspace{-2mm}
\begin{array}{lllll}
 \Kmin  - v^s              \\    
 \Kmax - v^s     \\
 \Kmin - v^s + y   \\ 
\Kmax - v^s + y 
\end{array}  \hspace{-1mm} \right ]\hspace{-1mm} .
%
%
%  = \left \{ \hspace{-2mm}
%\begin{array}{lllll}
%(0, \Kmin  - v^s )                    & \mbox{if }  i = 1 \\% \  b_1  > yp_{bs}   \\
%((\Kmin-v^s), (\Kmax - v^s))    & \mbox{if }   i = 2 \\% \  b_1 <yp_{bs} \le b_2 \\
%( ({\bar K}^s_Z  - v^s), (\Kmin - v^s + y))       \hspace{-5mm}& \mbox{if }  i = 3  \\%    b_2 < yp_{bs} \le   b_3   \\
%\Big ( \frac{ {\bar K}^s_Z   (1-p_{bs}) + (\Kmin +y) w p_{bs} }{1- p_{bs}(1-w) } - v^s,  %\hspace{-5mm}\\ \hspace{16mm} 
%(\Kmax - v^s + y) \Big)
%  \hspace{-1mm}& \mbox{if }  i = 4 \\% \   b_3  < yp_{bs} \le   b_4.
%\end{array} \right .
\end{eqnarray*}}
Let $b_5 :=y$, $c_5 := {\bar K}^s_Z  - v^s + y p_{bs}$ and $P_{D5} (p_{bs}) = 1$.
The common clearing  aggregate   (when {\small $b_{i-1} < yp_{bs} < b_{i} $}):
{\small $$
\bar{x_s}^\infty  =
   y (1-p_{bs}) -\frac{ \left ( yp_{bs} - c_i \right )  (1-  p_{bs})P_{Di}(p_{bs}) }{1 - (1-p_{bs})P_{Di}(p_{bs})} \ \forall i \le 5.   
$$}
The above is true when the big bank does not default, i.e., if  
{\small $$(k^b-\delta z_c - z_b)^+ - v^b + \bar{x_s}^\infty \frac{p_{bs}}{1-p_{bs}} > y p_{bs} $$}
\end{lemma}
{\bf Proof :} available in Appendix. \eop

  The above result indicates the `phase transitions' with respect to the connection parameter $ p_{bs}$. This lemma is true  
for the sub-case when big bank does not default (one example scenario,  when $(k^b - \delta z_c - z_b)^+ - v^b > y$).  However some of the `phase transition' results mentioned below are also true for the other case, by Lemma \ref{Lemma_First}.

\begin{figure} 
  \vspace{-17mm}
       \centering
       \hspace{-13mm}
       \begin{subfigure}
               \centering
               \includegraphics[width=4cm,height=6cm]{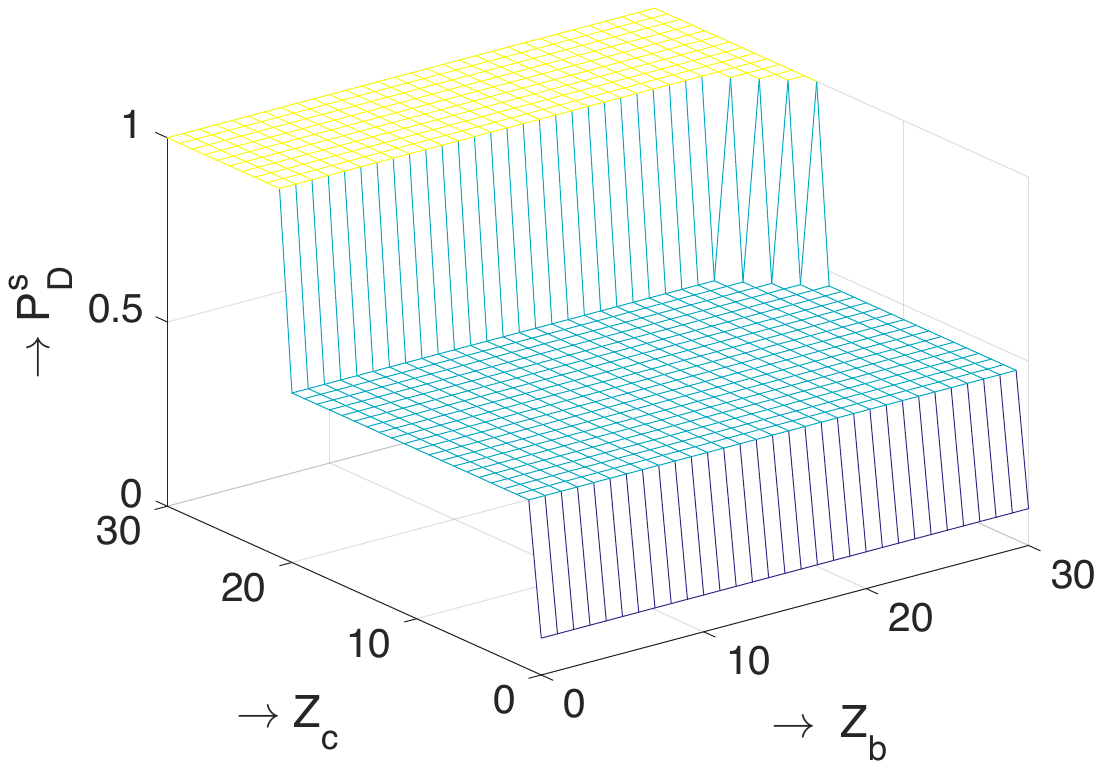}
               %\caption{Expected Surplus till time $T=1$}
       \end{subfigure} \hspace{-1mm}
       \begin{subfigure} 
               \centering
               \includegraphics[width=4cm,height=6cm]{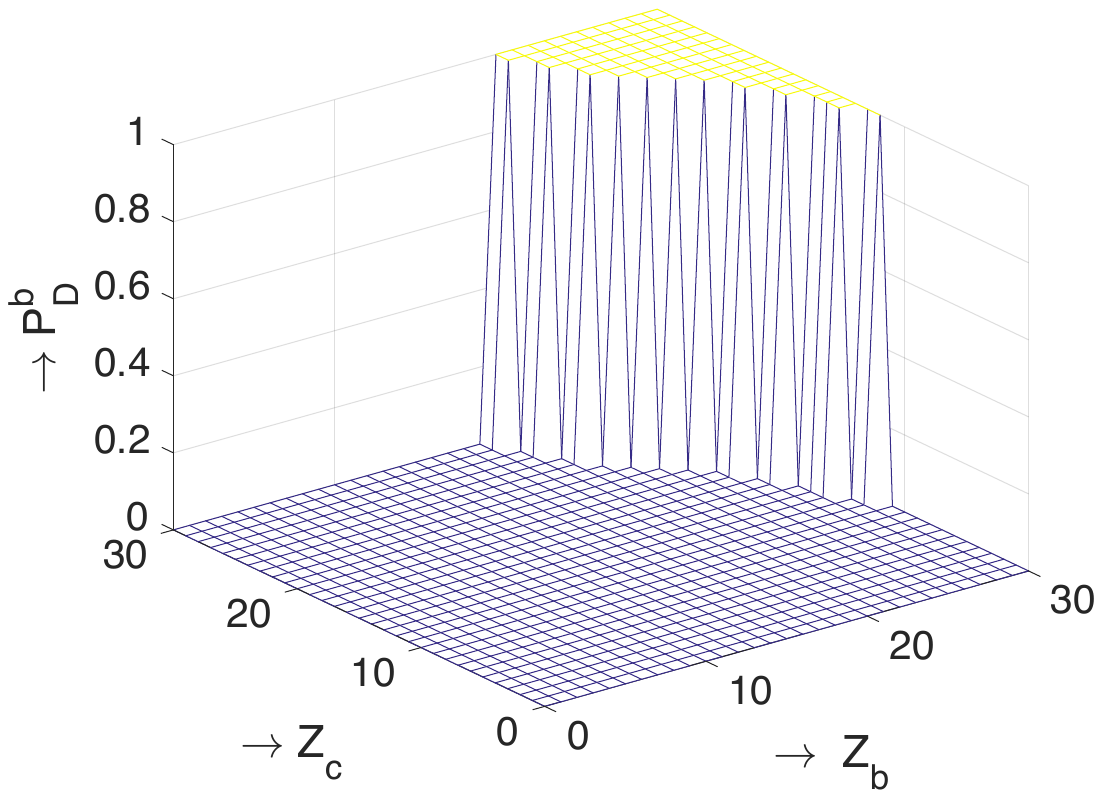}
            %  \caption{Expected fraction of defaults}
       \end{subfigure} \hspace{-15mm}
       \vspace{-22mm}
       \caption{For different shock realizations: {\small $ v^s \hspace{-1mm}= \hspace{-1mm}v^b \hspace{-1mm}=  \hspace{-1mm}12$,}  {\small $Z_i \hspace{-1mm}\sim \hspace{-1mm}Bin (.4,  20)$,     $(y, k^b, k^s) = (80, 55, 25)$,   $ (\delta, p_{bs}) = (.4 , .9)$}
   \label{Fig_ZbZc}}
  \vspace{-5mm}
\end{figure}

   As already discussed  when the connectivity parameter $p_{bs}$   is below $b_{1}/ y$, the network of all  small banks remains stable (the fraction of  defaults is zero).  
  But as soon as the connection parameter crosses  $b_{1}/ y = (\Kmin - v^s)/y$  the  small banks start defaulting, 
and we see a sharp jump of size   (see $P_{D1}$ and $P_{D2}$ in Lemma \ref{Lemma_BinaryNoRho}):
$$w(1-p_{bs}) =\frac{ w(y - \Kmin  + v^s)}{y} \mbox{  \underline{at  exactly} } p_{bs } = \frac{\Kmin - v^s}{ y}.$$
  These kind of phase-transitions can also be seen in Figures \ref{Fig_All_retuns}-\ref{Fig_ZbZc}.
  When  $p_{bs}$ is increased further, when it crosses the threshold $ b_{2}/y$,   $P_D^s$ has another jump/phase-transition.  
  At   $p_{bs} = b_2 / y$, we notice a sharp jump of size  (see $P_{D2}$ and $P_{D3}$ in Lemma \ref{Lemma_BinaryNoRho}):
  $$\left (1 - \frac{b_2}{y} \right ) (1-w). $$ 
 For this case, one needs to solve the equation  (Lemma \ref{Lemma_BinaryNoRho})
\begin{eqnarray*}
p_{bs} &=& \frac{b_2} { y} \\ 
& & \hspace{-10mm}=   \frac{(\Kmin-v^s) (1-p_{bs})^2w + (\Kmax -v^s) ( 1 - (1-p_{bs})^2w) }{y} 
\end{eqnarray*} to get the exact point of phase transition. 
 In a similar way from  Lemma \ref{Lemma_BinaryNoRho}, we see four possible phase transitions with respect to parameter $p_{bs}$. Since all the co-efficients also depend upon the shock 
 realizations ($z_c, z_b$) one can also obtain phase transitions with respect to shock sizes. Same is the case with other parameters.  
 
 More general observation from   Lemma \ref{Lemma_BinaryNoRho} is that, we have a possibility of small fraction of defaults when $p_{bs}$ is near 0 ($P_D^s $ can be 0)
or when $p_{bs}$ is near 1 with $b_4 > 1$ ($P_D^s \propto (1-p_{bs})$ or $\propto (1-p_{bs}(1-w))$).  The coefficient $b_4 $ with  $p_{bs}$ close to one, approximately 
equals $d_4 = \Kmax - v^s + y$, is mostly bigger than one.    Thus either small connectivity or large connectivity is better, but   intermediate connectivity may  not be good.  This can also be observed in the figures.

    \ForTR{
    \subsubsection*{With three time periods}
From  equation (\ref{Eqn_Finanace_FP}),  If $\Kmin +\rho^sA^s > v^s + y p_{bs}$  then ${\bar x}^{\infty*} = y (1- p_{bs})$   and 
 further $k^b > \delta z_c + z_b + v^b$ implies  $x_b^{\infty*} = y$. 
Thus the banks (small) may default, but they are able to clear their obligations completely by breaking the bonds.  One can easily verify the following in this case.
\begin{lemma}  
\label{Lemma_Binary2}
   Assume $k^b > \delta z_c + z_b + v^b$, $\epsilon < y$ and $\Kmin +\rho^sA^s > v^s + y p_{bs}$. Then  
\begin{eqnarray*}
P_D^s = \left \{ 
\begin{array}{lllll}
w (1-p_{bs})   & \mbox{if} &  \frac{ \Kmin  - v^s    }{y} < p_{bs} \le \frac{ \Kmax -   v^s   }{y}  \\
1-p_{bs}            & \mbox{if} &  \frac{ \Kmax -   v^s   }{y} < p_{bs} \le   \frac{ \Kmin + y - v^s    }{y}   \\
1- p_{bs}+w p_{bs}        & \mbox{if} &  \frac{ \Kmin + y - v^s    }{y} < p_{bs} \le   \frac{ \Kmax +y -   v^s   }{y}  \\
1                          & \mbox{if} &  \frac{ \Kmax +y -   v^s   }{y}  < p_{bs}  .
\end{array} \right .
\end{eqnarray*} 
\end{lemma}
  {\bf Proof:}
  When $y > \epsilon$ (i.e., $x_b > \epsilon$),  $ \Kmax < \Kmin + x_b $
   and 
  one can derive the result considering various scenarios as in the previous lemma.
  \eop\\
We continue with  sub-case considered of Lemma \ref{Lemma_Binary2} for which $P_D^b = 0.$
Define, 
\vspace{-4mm}
$$
P^{*s}_D  :=  \inf_{p_{bs} }  P^s_D, 
$$
the  minimum `expected defaults'  possible.  It is clear that for the sub-case considered in Lemma \ref{Lemma_Binary2}  it equals the following:
\begin{eqnarray*}
P^{*s}_D = \hspace{-8mm} &  \\  &  \left \{  \begin{array}{lllll}
%0   &   \mbox{ if }   \Kmin  > v^s \\
\min \left \{  w  \left (1- \frac{  \Kmax-v^s )} {y} \right),  \  \frac{ (v^s- \Kmin )^+} {y} \right \}  &\mbox{if } %\Kmin <
 v^s < \Kmax, \\
\min \left \{  w + (1-w) \frac{  v^s -\Kmax} {y},  \  \frac{ (v^s- \Kmin )} {y}, \ 1 \right \}  \hspace{-2mm} &\mbox{if }    v^s > \Kmax. 
\end{array} \right .
\end{eqnarray*}
Further  the assumptions of Lemma \ref{Lemma_Binary2} for any $p_{bs} = E[\eta^{bs}]$, we have  ${\bar x}^{\infty*} = y (1- p_{bs})$, $x_b^{\infty*} = y$,  and
 hence 
$$
  \Psi^s  = 
  (k^s -z_c- Z^s_{i})^++       y (    \eta^{bs}_i - E[\eta^{bs}]  ) -v^s  $$  and $
    \Psi^b  =  \psi^b =  k^b - z_c  \delta - z_b      -v^b > 0.  
$ 
Further $\Psi^s + \rho^sA^s \ge 0$ almost surely and 
thus  from equations (\ref{Eqn_ES1})-(\ref{Eqn_ES2}) the expected surplus equals
\begin{eqnarray*}
E[S(2)] &=&  E [ \Psi^s + \rho^s A^s  ; \Psi^s  <  0 + A^s ;  \Psi^s  \ge  0 ]  + \psi^b \\ 
& & \hspace{-20mm} \ =  \  E [   (k^s -z_c- Z^s_{i})^+   ] - v^s +  A^s  - (1-\rho^s)   A^s  P_D^s   + \psi^b.
\end{eqnarray*}
  Therefore the expected surplus is maximized at the same $p_{bs}^*$ which minimizes $P_D^s$  and then the optimal surplus is obtained by substituting  $P_D^s*$ in equation (\ref{Eqn_ES2}).

%  1-pbs - w (1-pbs) + w   = w  + (1-w) (1-pbs)

  }{}
  
\ignore { Check again
  When  $P_D^s = 1$ then 
  $$
  {\bar K}^s_Z + {\bar x} - v^s+  y p_{bs}   = {\bar x} / (1-p_{bs})
  $$
  Thus
    $$
 {\bar x} =  \frac{(y p_{bs} +  {\bar K}^s_Z  - v^s ) (1-p_{bs} ) } { 1 - (1-p_{bs}) P_D^s} $$

And hence
  $$
 {\bar x} = y (1-p_{bs})   + \frac{( {\bar K}^s_Z  - v^s ) (1-p_{bs} ) } { 1 - (1-p_{bs}) P_D^s} $$

 }

\underline{Influence of connectivity, shocks:} 
   We study the influence of connectivity parameter $p_{bs} = E[\eta^{bs}] = E[\eta^{sb}]$ for various shock scenarios. 
We begin with the case when $p_{bs} = 0^+$, i.e., as $p_{bs}$ approaches 0 (or is  $\approx 0$). There is negligible connection between the small banks and big bank, and    the small banks are primarily liable to other small banks.   The limit FP equations are: $ x_b = 0^+ $ and 

\vspace{-9mm}
{\small  \begin{eqnarray*}
  {\bar x}^{\infty*}  \hspace{-1.5mm} =\hspace{-.5mm} E \left [\hspace{-.5mm} \min \{   (k^s -z_c- Z^s_{i})^++    {\bar x}^{\infty*}  \hspace{-1mm} \ForTR{ +  \rho^s A^s}{} + \eta_i^{bs} 0^+ - v^s,  y \}  \right ]. 
  \end{eqnarray*}}
  When $\Kmin > v^s$,  it is clear that  $ {\bar x}^{\infty*} = y$, i.e.,  the small banks do not default at  $p_{bs} = 0^+$.  If  the big bank  encounters  big shock, i.e., if $( k^z - \delta z_c - z_b )^+ < v^b$, then 
  it  defaults near $p_{bs} = 0^+.$  In fact from equation (\ref{Eqn_Finanace_FP}),   the big bank defaults   for any $p_{bs} > 0$  under these conditions\footnote{
  Clearly {\small $ {\bar x}^{\infty*} \le  y (1-p_{bs})$},  {\small $y^b =  y (1-p_{bs}$)} and  so
  {\small $x_b < y$}.}. 
 Therefore we conclude,  {\it`when the big bank receives large shocks, it always defaults, the connection with the small banks does not help'}.   
Further  some of the small banks can also default leading to an increased fraction of defaults as $p_{bs}$ increases.   
  We considered one such example in the  first sub-figure of Figure \ref{Fig_All_retuns}, which reaffirms our observation: the fraction of defaults $P_D^s = 0$  for small $p_{bs}$ (till 0.1) 
  and near  $p_{bs} = 1$.  The  defaults  are more  for intermediate  $p_{bs}$.

 We consider the reverse situation now. Say the big bank does not default at $p_{bs} = 0^+$.  This is because it received small shocks such that  $( k^z - \delta z_c - z_b )^+ > v^b$.
  Say  some small banks receive big shocks such that
    $\Kmin < v^s$.  Then from  Lemma \ref{Lemma_BinaryNoRho}, $P_D^s \ge w$ at $p_{bs } = 0^+$.  
 The connection with big bank can improve the fraction of defaults, 
for example,  if we manage to chose a large enough $p_{bs}$ which  is between  $b_2/y$ and $b_3/y$ of  Lemma \ref{Lemma_BinaryNoRho} (and if further the big bank does not default).  In this case the asymptotic fraction of defaults could be smaller than $w$. 
We consider one such example   in the second sub-figure of Figure \ref{Fig_All_retuns}. We notice  that
$P_D^s$ reduces as $p_{bs}$ increases beyond 0.8, in fact 
 $P_D^s = 0$ (i.e, no small bank defaults) at  $p_{bs} = 1$.  
  Thus we conclude {\it `the big bank (with small shocks) can help the small banks'.}  
 
We consider a third example in Figure \ref{Fig_All_retuns}, where big bank and some small banks receive big shocks to default at $p_{bs} = 0^+$.
%Here we increased the individual shock of the big bank to 24.225.
The big bank continues to default at all $p_{bs}$, however  $P_D^s$ decrease  with increase in   $p_{bs}$. 
   Thus  for the chosen example of  economy with one big bank and many small banks, {\it `small banks can be stabilized by connecting to big bank, even if the later defaults, however they can't help the big bank'.}

\underline{Influence of shocks:} The banks can face   two types of shocks.   The idiosyncratic shocks ($\{Z_i^s\}, Z^b$)  are bank specific shocks,  while the common shock ($Z_c$)  affects all the  banks.  We aim to study  the role of magnitude of these shocks on the  cascading of defaults  in Figure \ref{Fig_ZbZc} for a fixed $p_{bs}=0.9$. We observe  two phase transitions in the fraction of defaults, $P_D^s$ remains at 0.1  ($1-p_{bs}$) for some region of $(Z_b, Z_c)$, jumps to 
0.46 ($1-p_{bs} (1-w)$) 
 and then to 1, and, one phase transition with respect to big bank.   This behaviour can also be explained using the barrier constants $\{b_i \}$ of Lemma \ref{Lemma_BinaryNoRho}, which depend   upon these shocks.  One can observe significant sharp jumps at the phase transition points,  and these points are very important for any financial network. These jumps and transitions points can be studied either using Lemma \ref{Lemma_BinaryNoRho} or by studying  the   simplified FP equations (\ref{Eqn_Finanace_FP}) numerically. 

\
\section*{Conclusions}

 We considered a  random graph, with edges representing the influence factors between a big  (highly influential)  node and numerous  small nodes.   
 The performance/status of individual nodes is resultant of these  influences, which are represented by  fixed point (FP) equations. We   showed that the solution of the random FP equations  converge almost surely to that of a limit system and these solutions are asymptotically independent. One may have  multiple solutions 
for finite graphs, however any sequence of them  converge to the unique FP of the limit system (if it has unique FP).   Thus we have a procedure to solve the large dimensional FP equations, using mean-field kind of techniques.  The proposed solution requires solving of  'aggregate' FPs in a much smaller dimensional space and is accurate asymptotically.

The clearing vectors (the fraction of liabilities eventually cleared) in a financial network are generally represented by  random FP equations and we studied the same using our results.
We study  an example  heterogeneous  financial network with one big bank and many small banks.
We have reduced the overall economy problem in this set-up to a two node problem - one big bank and one aggregate small bank, thus facilitating big picture analysis.
We observe some interesting phase transitions, one can easily study the nature of these phase transitions using the approximate solutions of the involved FPs.  
%We also observe that 
%the small banks can be stabilized by connecting with big banks and not the vice versa.  
When small banks invest more in big banks, lesser fraction of them default and this is true even when 
all of them face  large idiosyncratic shocks.
 These observations could be specific to the example considered by us,  however we now have a procedure to study complex networks and a more elaborate  study would help us derive more concrete observations.   
One can easily generalize the results by relaxing many of the assumptions, one can   apply this approach to   more applications and these two would be  the topics of future interest.

%As already mentioned, one can think of $X_b, \{ X_i \}_i$ as some performance measures that are influenced by the performance of other nodes and the above set of fixed point equations are valid when the influence  of small nodes is via a random weighted  aggregate  of  their performances.  Our aim is to solve these fixed point equations, in an asymptotic limit obtained by increasing the number of nodes $n$ to infinity. As of now these equations consider small nodes that are similar and one b-node. One can easily generalize to the case when there are finite 
%number of groups of small nodes, wherein the nodes within a group are identical (in  distributional sense).

\section{Markov decision process (MDP) - state aggregation}
Say we have finite number of actions in ${\cal A}$. There exists finite number of groups $\{G_l \}_{l \le L}$
and
if the states of an infinite MDP aggregate in the following manner
$$
\bigg|\sum_{j \in G_l} p^{(n)}( j | i, a) - p( l |k,  a)  \bigg|  \to 0\mbox{, i.e., } $$
$$ \bigg|  p^{(n)}( G_l | i, a) - p( l |k,  a)  \bigg|  \mbox{ for all }  l \in G_k
$$
and if  the immediate rewards also converge
$$
r^{(n)}(i, a) =  r(k, a)  \mbox{ for all }  i \in G_k.
$$
Then using our theorem we can show that the value functions 
$$
v^{(n)} (i) = \min_{a \in {\cal A} } \{r(i, a) +  \lambda \sum_{j} p^{(n)} (j| i, a) v^{(n)} ( j) \}
$$
converges
$$
v^{(n)} (i) \to v(k)  \mbox{ for all }  i \in G_k,
$$
and also the optimal strategy converges 
$$
a^{(n)*} (i) \to a^{*} (k) .
$$
This will be true if the limit system has unique optimizer.

Idea is to use convergence of aggregates
$$
{\bar v}^{(n)}_{i, a} :=  \sum_{j} p^{(n)} (j| i, a) v^{(n)} ( j)  =   \sum_{j} p^{(n)} (j| i, a)  \xi_j  (  {\bf {\bar v}}^{(n)}_j )
$$ where  the vector of aggregate for any $i$ is defined as:
$$
{\bf {\bar v}}^{(n)}_{i} := \{ {\bar v}^{(n)}_{i,a} \}_a,
$$
and then for all $j \in G_l$
$$
 \xi_j  (  {\bf {\bar v}}^{(n)}_j ) := \min_a  \left \{ r (l, a) +  {\bar v}^{(n)}_{j, a} ) \right   \}.
$$
\section*{Appendix: Proof of Lemma \ref{Lemma_BinaryNoRho} }
\textbf{Proof of Lemma \ref{Lemma_BinaryNoRho}:} The big bank does not default hence $ { x}^{\infty *}_b=y$.
 Therefore we  can rewrite the  FP equation representing aggregate clearing vector $\bar{x_s}^\infty$  as  below:
\begin{eqnarray}
      \frac{\bar{x_s}^\infty}{1-p_{bs}} &=& 
      \min\lbrace\Kmin -v^s+ \bar{x_s}^\infty,y\rbrace w(1-p_{bs}) + \min\lbrace\Kmax -v^s+ \bar{x_s}^\infty,y\rbrace (1-w)(1-p_{bs}) \nonumber
\\
    &&   + \min\lbrace\Kmin -v^s +\bar{x_s}^\infty + y,y\rbrace wp_{bs}+ 
        \min\lbrace\Kmax -v^s +\bar{x_s}^\infty+y,y\rbrace (1-w)p_{bs}. \label{Eqn_FP_Binary}
      \end{eqnarray}

It is clear that $\Kmin < \Kmax$ and $\Kmin < \Kmin + y$ etc. 
If $ y > \epsilon$ then we also have $\Kmax <  \Kmin + y$. Then the above FP equation has a natural order in the following sense: the terms are arranged in increasing order when the corresponding probabilities are not considered. For example  the third term,  $\min\lbrace\Kmin -v^s +\bar{x_s}^\infty + y,y\rbrace  \le \min\lbrace\Kmax -v^s +\bar{x_s}^\infty+y,y\rbrace$, the fourth term.
The best scenario is with fourth term  (small  banks receive   zero shock and connect with big bank) while the worst is with the first term (small bank face negative shock and are  not connected to   big bank).

 \textbf{Case 1:}  There is no default even in the worst   scenario i.e. if
\begin{eqnarray*}
  \Kmin -v^s+ \bar{x_s}^\infty > y .  
 \end{eqnarray*}
 Then none of the  small banks  default leading to $P_D^s = P_{D1} =0$,  and  hence
  the clearing vector satisfies $
\frac{\bar{x_s}^\infty}{1-p_{bs}} = y,
 $
 or equivalently $\bar{x_s}^\infty =y(1-p_{bs})$. That is  Case 1 holds as long as:
 \begin{eqnarray}
  \Kmin -v^s+ y (1- p_{bs}) > y  \mbox{ or equivalently as long as }  y p_{bs } < \Kmin -v^s. \label{Eqn_case1}
 \end{eqnarray}
  However as $p_{bs}$ increases, the above may not be true and this gives us the bound   $b_1 = \Kmin -v^s.$\\
% \begin{align*}
%  \Kmin -v^s+ \bar{x_s}^\infty = y 
%  \Rightarrow \Kmin -v^s+ y(1-p_{bs})=y
%   \Rightarrow yp_{bs}= \Kmin -v^s
%    \end{align*}
  \textbf{Case 2}  When there is default only in the first term of (\ref{Eqn_FP_Binary}),  i.e., when $P_D^s (p_{bs}) = P_{D2} (p_{bs})=w(1-p_{bs})$.
  The aggregate clearing vector in this case satisfies:
      \begin{align*}
      \frac{\bar{x_s}^\infty}{1-p_{bs}} = 
      (\Kmin -v^s+ \bar{x_s}^\infty)P_{D2}+ y(1-P_{D2})
    \\  \Rightarrow {\small 
\bar{x_s}^\infty  =
   y (1-p_{bs}) -\frac{ \left ( yp_{bs} - c_2 \right )  (1-  p_{bs})P_{D2}(p_{bs}) }{1 - (1-p_{bs})P_{D2}(p_{bs})}  \mbox{ where } c_2= \Kmin -v^s.   
}
      \end{align*}
   The Case 2 holds as long as   
 \begin{eqnarray}
  \Kmin -v^s+ \bar{x_s}^\infty  <  y  \mbox{ and  }   \Kmax -v^s+  \bar{x_s}^\infty  >  y.  
 \end{eqnarray}
 Once again as $p_{bs}$ increases, the above may not be true (the second inequality  can fail) and this gives us the bound    $b_2$.    
  The bound $b_2$ can be obtained:
 \begin{align*}
  \Kmax -v^s+ \bar{x_s}^\infty = y 
   \Rightarrow yp_{bs}=(\Kmin-v^s)P_{D2}(1-p_{bs}) +(\Kmax -v^s)(1-(1-p_{bs})P_{D2})
 \end{align*}
   \begin{align*}
\Rightarrow yp_{bs}= c_2 P_{D2}(1-p_{bs})+d_2(1-(1-p_{bs})P_{D2}) \mbox{ where } d_2= \Kmax -v^s.
    \end{align*}
 
     Thus bound, $b_2 = c_2 P_{D2}(1-p_{bs})+d_2(1-(1-p_{bs})P_{D2}).$

     \textbf{Case 3}  When there is default only in the first two  terms  of (\ref{Eqn_FP_Binary}),  i.e., when $P_D^s (p_{bs}) = P_{D3} (p_{bs})=(1-p_{bs})$.
  The aggregate clearing vector in this case satisfies:
      \begin{align*}
      \frac{\bar{x_s}^\infty}{1-p_{bs}} = 
      (\Kmin -v^s+ \bar{x_s}^\infty)w(1-p_{bs})+   (\Kmax -v^s+ \bar{x_s}^\infty)(1-w)(1-p_{bs})+ y(1-P_{D3}).
      \end{align*}
      This implies 
     \begin{align}
     \label{Eqn_XsCase3}
\bar{x_s}^\infty  =
   y (1-p_{bs}) -\frac{ \left ( yp_{bs} - c_3 \right )  (1-  p_{bs})P_{D3}(p_{bs}) }{1 - (1-p_{bs})P_{D3}(p_{bs})}  \mbox{ where } c_3= {\bar K}^s_Z  - v^s.   
      \end{align}
   The Case 3 holds as long as   
 \begin{eqnarray}
  \Kmax -v^s+  \bar{x_s}^\infty  <  y  \mbox{ and  }   \Kmin -v^s+ \bar{x_s}^\infty +y >  y.  
 \end{eqnarray}
 Using (\ref{Eqn_XsCase3}) one can easily show that that $\bar{x_s}^\infty$ is decreasing with increase in $p_{bs}$.
Thus again as $p_{bs}$ increases, the above inequalities (second one) may not be true and this gives us the bound    $b_3$.    
  The bound $b_3$ can be obtained:
 \begin{align*}
  \Kmin -v^s+ \bar{x_s}^\infty+y = y 
   \Rightarrow yp_{bs}= c_3 P_{D3}(1-p_{bs})+d_3(1-(1-p_{bs})P_{D3})
 \end{align*}
   \begin{align*}
\Rightarrow yp_{bs}= c_3 P_{D3}(1-p_{bs})+d_3(1-(1-p_{bs})P_{D3}) \mbox{ where } d_3= \Kmin -v^s+y.
    \end{align*}
 
     Thus bound, $b_3 = c_3 P_{D3}(1-p_{bs})+d_3(1-(1-p_{bs})P_{D3}).$
     
     Continuing this way one can obtain all the sub-cases of the lemma.

     \eop

\end{document}